\documentclass{article}[12pt]
\usepackage{amsmath,amssymb}
\usepackage{amscd}
\usepackage[all]{xy}

\newtheorem{thm}{Theorem}
\newtheorem{df}[thm]{Def\/inition}
\newtheorem{prop}[thm]{Proposition}
\newtheorem{cor}[thm]{Corollary}
\newtheorem{lemma}[thm]{Lemma}

\newcommand{\eqn}{\begin{equation}}
\newcommand{\eeqn}{\end{equation}}
\newcommand{\ideal}{\mathfrak}

\DeclareMathOperator{\ord}{ord} 
 
\DeclareMathOperator{\Hom}{Hom} 
\DeclareMathOperator{\Ext}{Ext} 
\DeclareMathOperator{\id}{id}
\DeclareMathOperator{\gr}{gr}  
\DeclareMathOperator{\Spec}{Spec}
\DeclareMathOperator{\Spf}{Spf}

\DeclareMathOperator{\perf}{perf}  
\DeclareMathOperator{\length}{length}  
\DeclareMathOperator{\depth}{depth}  

\def\endproof{$\hfill \square$}

\begin{document}
\title{Purity results for $p$-divisible groups and abelian schemes over regular bases of mixed characteristic}
\author{Adrian Vasiu and Thomas Zink}
\maketitle
\centerline{\it To Michel Raynaud, for his 71th birthday. To appear in Doc. Math.}
\bigskip

\medskip\noindent
{\bf Abstract.} Let  $p$ be a prime. Let $(R,\ideal{m})$ be a regular local ring of mixed
characteristic $(0,p)$ and absolute index of ramification $e$. We provide general
criteria of when each abelian scheme over $\Spec R\setminus\{\ideal{m}\}$ extends to an abelian scheme over $\Spec R$. We show that such extensions always exist if $e\le p-1$, exist in most cases if $p\le e\le 2p-3$, and do not exist in general  if $e\ge
2p-2$. The case $e\le p-1$ implies the uniqueness of
integral canonical models of Shimura varieties over a discrete
valuation ring $O$ of mixed characteristic $(0,p)$ and index of
ramification at most $p-1$. This leads to large classes of examples of
N\'eron models over $O$. If $p>2$ and index $p-1$, the examples are
new. 

\bigskip\noindent
{\bf Key words}: rings, group schemes, $p$-divisible groups, Breuil windows and modules, abelian schemes, Shimura varieties, and N\'eron models.

\bigskip\noindent
{\bf MSC 2000:} 11G10, 11G18, 14F30, 14G35, 14G40, 14K10, 14K15, 14L05, 14L15, and
14J20.

\section{Introduction}
Let $p$ be a prime number. We recall the following {\it global purity} notion introduced in \cite{V1}, Definitions 3.2.1  2) and 9) and studied in \cite{V1} and
\cite{V2}. 

\begin{df}\label{D1}
Let $X$ be a regular scheme that is faithfully flat over $\Spec \mathbb{Z}_{(p)}$. 
We say $X$ is healthy regular (resp. $p$-healthy regular), if for each
open subscheme $U$ of $X$ which contains $X_{\mathbb{Q}}$ and all
generic points of $X_{\mathbb{F}_p}$, every abelian scheme
(resp. $p$-divisible group) over $U$ extends uniquely to an abelian scheme (resp. a
$p$-divisible group) over $X$. 
\end{df}

In (the proofs of) \cite{FC}, Chapter IV, Theorems $6.4$, $6.4'$, and
$6.8$ was claimed that every regular scheme which is faithfully flat
over $\Spec \mathbb{Z}_{(p)}$ is healthy regular as well as $p$-healthy
regular. This claim was disproved by an example
of Raynaud--Gabber (see \cite{Ga} and \cite{dJO}, Section 6): the regular scheme
$\Spec W(k)[[T_1,T_2]]/(p-(T_1T_2)^{p-1})$
 is neither $p$-healthy nor healthy regular.  Here $W(k)$ is the ring
 of Witt vectors with coefficients in a perfect field $k$ of
 characteristic $p$. The importance of healthy and $p$-healthy regular
 schemes stems from their applications to the study of integral models
 of {\it Shimura varieties}. We have a {\it local} version of
 Definition \ref{D1} as suggested by Grothendieck's work on the classical
 Nagata--Zariski purity theorem (see \cite{Gr}).  

\begin{df}\label{D2}
Let $R$ be a local noetherian ring with maximal ideal $\ideal{m}$
such that $\depth R \geq 2$. We say that $R$ is quasi-healthy (resp.
$p$-quasi-healthy) if each abelian scheme (resp. $p$-divisible group) over  $\Spec R
\setminus \{\ideal{m}\}$ extends uniquely to an abelian scheme (resp. a
$p$-divisible group) over $\Spec R$. 
\end{df}

\noindent
If $R$ is local, complete, regular of dimension $2$ and mixed characteristic $(0,p)$, then the fact that $R$ is $p$-quasi-healthy can be restated in terms of finite flat commutative group schemes annihilated by $p$ over $\Spec R$ (cf. Lemma \ref{Le3}). 

Our main result is the following theorem proved in Subsections  \ref{S3} and
\ref{S5}.

\begin{thm}\label{T1}
Let $R$ be a regular local ring of dimension $d\geq 2$ and of
mixed characteristic $(0,p)$. We assume that there exists a faithfully flat
local $R$-algebra $\hat{R}$ which is complete and regular of dimension $d$, which
has an algebraically closed
residue class field $k$, and which is equipped with an epimorphism $\hat{R}
\twoheadrightarrow W(k)[[T_1,T_2]]/(p - h)$ where $h \in (T_1,
T_2)W(k)[[T_1,T_2]] $ is a power series whose reduction modulo the
ideal $(p,T_1^p,T_2^p,T_1^{p-1}T_2^{p-1})$ is non-zero. Then  $R$ is 
quasi-healthy. If moreover $d=2$, then $R$ is also $p$-quasi-healthy.
\end{thm}
For instance, Theorem \ref{T1} applies if the strict completion of $R$ is isomorphic
to $W(k)[[T_1,\ldots,T_d]]/(p-T_1\cdot\ldots\cdot T_m)$ with $1\leq
m\le\min\{d,2p-3\}$ (cf. Subsection \ref{S6}). The following consequence is also
proved in Subsection \ref{S5}.

\begin{cor}\label{T2}
Let $R$ be a regular local ring of dimension $d \geq 2$ and of
mixed characteristic $(0,p)$. Let $\ideal{m}$ be the maximal ideal of $R$. We
assume that $p \notin \ideal{m}^p$. Then $R$ is quasi-healthy. If moreover $d =
2$, then $R$ is also $p$-quasi-healthy.
\end{cor}

Directly from Theorem \ref{T1} and from very definitions we get:

\begin{cor}\label{C1} 
Let $X$ be a regular scheme that is faithfully flat over
$\Spec \mathbb{Z}_{(p)}$. We assume that each local ring $R$ of $X$ of mixed
characteristic $(0,p)$ and dimension at least $2$ is such that the hypotheses of
Theorem \ref{T1} hold for it (for instance, this
holds if $X$ is formally smooth over the spectrum of a discrete valuation ring $O$ of
mixed characteristic $(0,p)$ and index of ramification $e\leq p-1$). Then $X$ is
healthy regular. If moreover $\dim X=2$, then $X$ is
also $p$-healthy regular. 
\end{cor}

The importance of Corollary  \ref{C1} stems from its applications to {\it
N\'eron models} (see Section 6). Theorem  \ref{T5} shows the existence of large
classes of new types of N\'eron models that were not studied before in \cite{N},
\cite{BLR}, \cite{V1}, \cite{V2}, or \cite{V3}, Proposition 4.4.1. Corollary
\ref{C1} encompasses (the correct parts of) \cite{V1}, Subsubsection 3.2.17 and
\cite{V2}, Theorem 1.3.
Theorem \ref{T4} (i) shows that if $X$ is formally smooth over the spectrum of a
discrete valuation ring $O$ of mixed characteristic $(0,p)$ and index of
ramification $e$ at least $p$, then in general $X$ is neither $p$-healthy nor
healthy regular. From this and Raynaud--Gabber example we get that Theorem \ref{T1}
and Corollary \ref{T2} are optimal. Even more, if $R=W(k)[[T_1,T_2]]/(p-h)$ with
$h\in (T_1,T_2)$, then one would be inclined to expect that $R$ is $p$-quasi-healthy if
and only if $h$ does not belong to the ideal $(p,T_1^p,T_2^p,T_1^{p-1}T_2^{p-1})$; this is supported by Theorems \ref{T1} and \ref{T4} and by Lemma
\ref{Le2}. In particular, parts (ii) and (iii) of Theorem \ref{T4} present two
generalizations of the Raynaud--Gabber example.

Our proofs are based on the classification of finite flat commutative group schemes of $p$ power order over the spectrum of a local, complete, regular ring $R$ of mixed characteristic $(0,p)$ and perfect residue class field. For $\dim R = 1$ this classification was a conjecture of Breuil \cite{Br} proved by Kisin in \cite{K1} and \cite{K2} and reproved by us in \cite{VZ}, Theorem 1. Some cases with $\dim R \geq 2$ were also treated in \cite{VZ}. The general case is proved by Lau in \cite{L1},
Theorems 1.2 and 10.7. 
Proposition \ref{P4} provides a new proof of Raynaud's
result \cite{R2}, Corollary 3.3.6.

Subsection 5.1 disproves an additional claim of \cite{FC}, Chapter V, Section 6. It
is the claim of \cite{FC}, top of p. 184 on torsors of liftings of $p$-divisible
groups which was not previously disproved and which unfortunately was used in \cite{V1}
and \cite{V2}, Subsection 4.3. This explains why our results on $p$-healthy regular
schemes and $p$-quasi-healthy regular local rings work only for dimension $2$ (the
difficulty is for the passage from dimension $2$ to dimension $3$).
Implicitly, the $p$-healthy part of \cite{V2}, Theorem 1.3 is proved correctly in
\cite{V2} only for dimension $2$.

The paper is structured
as follows. Different preliminaries on {\it Breuil windows} and {\it modules} are
introduced in Section 2. In Section 3 we study morphisms between Breuil modules.
Our basic results on extending properties of finite flat group schemes,
$p$-divisible groups, and abelian schemes are presented in Sections 4 and 5. Section
6 contains applications to integral models and N\'eron models. 

\medskip
{\bf Acknowledgement.}
The first author would like to thank Binghamton and Bielefeld Universities for good
conditions with which to write this note; he was partially supported by the NSF grant DMS \#0900967. The second author would like to thank Eike
Lau for helpful discussions. Both authors would like to thank the referee for several valuable comments.

\section{Preliminaries}\label{S0}
In this paper the notions of {\it frame} and {\it window} are in a more
general sense than in \cite{Z1}. The new notions are suggested by the 
works \cite{Br}, \cite{K1}, \cite{VZ}, and \cite{L1}. In all that follows we
assume that a ring is unitary and commutative and that a finite flat group scheme is
commutative and (locally) of $p$ power order. 

\begin{df}\label{D3} A frame $\mathcal{F} = (R,S,J, \sigma, \dot{\sigma},
 \theta)$ for a ring $R$ consists of the following data: 
\begin{itemize}
\item[(a)] A ring $S$ and an ideal $J \subset S$.
\item[(b)] An isomorphism of rings  $S/J \cong R$.
\item[(c)] A ring homomorphism $\sigma: S \rightarrow S$.
\item[(d)] A $\sigma$-linear map $\dot{\sigma} : J
  \rightarrow S$.
\item[(e)] An element $\theta \in S$.
\end{itemize}
We assume that $pS +J$ is in the radical of $S$, that $\sigma$ induces
the Frobenius endomorphism on $S/pS$, and that the following equation holds:
\begin{equation}\label{frame1e}
   \sigma (\eta) = \theta \dot{\sigma} (\eta), \quad \text{for}\;\text{all}\; \eta
   \in J.
\end{equation}
\end{df} 

Assume that 
\begin{equation}\label{frame2e}
   S\dot{\sigma} (J) = S.
\end{equation}
 In this case we find an
equation $1 = \sum_i\xi_i \dot{\sigma}(\eta_i)$ with $\xi_i \in S$ and
$\eta_i \in J$. From this and (\ref{frame1e}) we get that $\theta = \sum_i \xi_i
\sigma (\eta_i)$. Therefore for each $\eta \in J$ we have 
$$\sigma(\eta) = \sum_i\xi_i \sigma(\eta) \dot{\sigma}(\eta_i) =  
\sum_i\xi_i \dot{\sigma}(\eta \eta_i) = \sum_i \xi_i \sigma(\eta_i)
\dot{\sigma}(\eta) = \theta \dot{\sigma}(\eta).$$
We conclude that the equation (\ref{frame2e}) implies the existence and uniqueness of an element $\theta$ such that the equation $(\ref{frame1e})$ is satisfied. 

If $M$ is an $S$-module we set $M^{(\sigma)} := S\otimes_{\sigma,
  S}M$. The linearization of a $\sigma$-linear map $\phi: M \rightarrow
N$ is denoted by $\phi^{\sharp}: M^{(\sigma)} \rightarrow N$.

\begin{df}\label{D4}
A window with respect to $\mathcal{F}$ is a quadruple
$(P,Q,F,\dot{F})$ where:
\begin{itemize}
\item[(a)] $P$ is a finitely generated projective $S$-module.
\item[(b)] $Q \subset P$ is an $S$-submodule. 
\item[(c)] $F : P \rightarrow P$ is a $\sigma$-linear map.
\item[(d)] $\dot{F}: Q \rightarrow P$ is a $\sigma$-linear map. 
\end{itemize}

We assume that the following three conditions are satisfied: 
\begin{itemize}
\item[(i)] There exists a decomposition $P = T \oplus L$ such that $Q = JT
  \oplus L$. 
\item[(ii)] $F(y) = \theta \dot{F}(y)$ for $y \in Q$ and $\dot{F}(\eta
  x) = \dot{\sigma}(\eta) F(x)$ for $x \in P$ and $\eta \in J$.
\item[(iii)] $F(P)$ and $\dot{F}(Q)$ generate $P$ as an $S$-module.
\end{itemize}
\end{df}

\indent
If (\ref{frame2e}) holds, then from the second part of $(ii)$ we get that for $x\in P$ 
we have $F(x)=1F(x)=\sum_{i} \xi_i\dot{\sigma}(\eta_i)F(x)=\sum_{i} \xi_i \dot{F}(\eta_i x)$; 
if moreover $x\in Q$, then $F(x)=\sum_{i} \xi_i \sigma(\eta_i)\dot{F}(x)=\theta \dot{F}(x)$. 
Thus if (\ref{frame2e}) holds, then we have an inclusion 
$F(P) \subset S\dot{F}(Q)$ and the first condition of $(ii)$ follows
from the second condition of $(ii)$. A decomposition as in $(i)$ is called
a {\it normal decomposition}.

We note that for each window $(P,Q,F,\dot{F})$ as above, the $S$-linear map 
\begin{equation}\label{windows1e}
   F^{\sharp} \oplus \dot{F}^{\sharp}: S \otimes_{\sigma,S} T \oplus
   S \otimes_{\sigma,S} L \rightarrow T \oplus L
\end{equation}
is an isomorphism.  Conversely an arbitrary isomorphism
(\ref{windows1e}) defines uniquely a window with respect to $\mathcal{F}$ equipped
with a given normal decomposition.

Let $R$ be a regular local ring of mixed characteristic $(0,p)$. Let
$\ideal{m}$ be the maximal ideal of $R$ and let $k := R/\ideal{m}$. 

A finite flat group scheme $H$ over $\Spec R$ is called
{\it residually connected} if its special fibre $H_k$ over $\Spec k$ is connected
(i.e.,
has a trivial \'etale part). If $R$ is complete, then $H$ is residually connected if
and only if $H$ as a scheme is connected and therefore in this case we will drop the
word residually. We apply the same terminology to $p$-divisible groups over $\Spec
R$.

In this section we will assume moreover that $k$ is perfect and that $R$ is
complete (in the $\ideal{m}$-adic topology). Let $d=\dim R\ge 1$.

We choose regular parameters $t_1, \ldots, t_d$ of $R$. We denote by
$W(\dag)$  the ring of Witt vectors with coefficients in a ring
$\dag$. We set $\mathfrak{S} := W(k)[[T_1, \ldots, T_d]]$. We consider
the epimorphism of rings
\begin{displaymath}
   \mathfrak{S} \twoheadrightarrow R,
\end{displaymath} 
which maps each indeterminate $T_i$ to $t_i$.

Let $h \in \mathfrak{S}$ be a power series without constant
term such that we have
\begin{displaymath}
   p = h(t_1, \ldots ,t_d).
\end{displaymath}
We set $E := p - h \in \mathfrak{S}$. Then we have a canonical isomorphism
\begin{displaymath}
   \mathfrak{S}/E\mathfrak{S} \cong R.
\end{displaymath}
As $R$ is a regular local ring of mixed characteristic $(0,p)$, we have $E\notin p\mathfrak{S}$. We extend the Frobenius automorphism $\sigma$ of $W(k)$ to $\mathfrak{S}$ by the rule
\begin{displaymath}
   \sigma (T_i) = T_i^p.
\end{displaymath}
Note that $\sigma:\mathfrak{S}\to\mathfrak{S}$ is flat. Let  $\dot{\sigma}:E\mathfrak{S}\rightarrow \mathfrak{S}$ be the $\sigma$-linear map
defined
by the rule 
\begin{displaymath}
   \dot{\sigma}(Es) := \sigma(s).
\end{displaymath} 
\begin{df}\label{D5}
We refer to the sextuple $\mathcal{L} = (R, \mathfrak{S}, E\mathfrak{S}, \sigma,
\dot{\sigma}, \sigma(E))$ as a standard frame for the ring $R$. A window for this
frame will be called a {\it Breuil window}. 
\end{df}

As $\dot{\sigma}(E)=1$, the condition (\ref{frame2e}) holds for $\mathcal{L}$. 
Let
$\hat{W}(R) \subset W(R)$ be the subring defined in \cite{Z2}, Introduction and let
$\sigma_R$ be
its Frobenius endomorphism. There exists a natural ring
homomorphism
\begin{displaymath}
   \varkappa : \mathfrak{S} \rightarrow \hat{W}(R)
\end{displaymath}
which commutes with the Frobenius endomorphisms, i.e. for all $s\in \mathfrak{S}$ we
have
\begin{displaymath}
   \varkappa (\sigma (s)) = \sigma_R(\varkappa (s)).
\end{displaymath}
The element $\varkappa (T_i)$ is the Teichm\"uller
representative $[t_i]\in \hat{W}(R)$ of $t_i$.

\begin{thm}\label{T3} (Lau \cite{L1}, Theorem 1.1)
If $p \geq 3$, then the category of Breuil windows is equivalent to the
category of $p$-divisible groups over $\Spec R$.
\end{thm}
This theorem was proved first in some cases in \cite{VZ}, Theorem
1. As in \cite{VZ}, the proof in  \cite{L1} shows first that the
category of Breuil windows is equivalent (via $\varkappa$) to the
category of Dieudonn\'e displays over $\hat{W}(R)$. Theorem \ref{T3}
follows from this and the fact (see \cite{Z2}, Theorem) that the category of
Dieudonn\'e displays over $\hat{W}(R)$ is equivalent to the
category of $p$-divisible groups over $\Spec R$.  

There exists a version of Theorem \ref{T3} for connected
$p$-divisible groups over $\Spec R$ as described in the introductions of
\cite{VZ} and \cite{L1}.  This version holds as well for $p = 2$, cf. \cite{L1}, Corollary 10.6.

The categories of Theorem \ref{T3} have natural exact structures \cite{Me2}. The
equivalence of Theorem \ref{T3} (or of its version for $p\geq 2$) respects the exact structures.

For a prime ideal $\ideal{p}$ of $R$ which contains $p$, let $\kappa
(\ideal{p})^{\perf}$ be the perfect hull of the residue class field $\kappa
(\ideal{p})$ of $\ideal{p}$. We deduce from $\varkappa$ a ring homomorphism
$\varkappa_{\ideal{p}}:\mathfrak{S}
\rightarrow W(\kappa(\ideal{p})^{\perf})=\hat{W}(\kappa(\ideal{p})^{\perf})$.

\begin{prop}\label{P1}
Let $G$ be a $p$-divisible group over $R$. If $p=2$, we assume that $G$ is connected. Let $(P,Q,F,\dot{F})$ be the Breuil window of $G$. Then the classical
Dieudonn\'e module of $G_{\kappa (\ideal{p})^{\perf}}$ is canonically
isomorphic (in a functorial way) to $W(\kappa (\ideal{p})^{\perf})
\otimes_{\varkappa_{\ideal{p}},\mathfrak{S}}
P$ endowed with the $\sigma_{\kappa(\ideal{p})^{\perf}}$-linear map $\sigma_{\kappa(\ideal{p})^{\perf}} \otimes F$.
\end{prop}
{\bf Proof.} 
If $G$ is connected, we consider its display $(P',Q',F',\dot{F}')$ over
$W(R)$. By \cite{L1}, Theorem 1.1 and Corollary 10.6 we have 
$P' = W(R) \otimes_{\varkappa,\mathfrak{S}} P$. Let us denote by
$\mathbb{D}(G)$ the Grothendieck--Messing crystal associated
to $G$, cf. \cite{Me1}. By \cite{Z1}, Theorem 6 there is a canonical and functorial isomorphism
$\mathbb{D}(G)_{W(R)} \cong  P' \cong W(R) \otimes_{\varkappa, \mathfrak{S}} P$. The functor $\mathbb{D}$ commutes
with base change. If we apply this to the morphism of pd-extensions
$W(R) \rightarrow  W(\kappa(\ideal{p})^{\perf})$ of $R$ and
$\kappa(\ideal{p})^{\perf}$ (respectively), we get that 
$\mathbb{D}(G_{\kappa (\ideal{p})^{\perf}}) \cong \mathbb{D}(G)_{W(\kappa(\ideal{p})^{\perf})}\cong W(\kappa(\ideal{p})^{\perf})
\otimes_{\varkappa_{\ideal{p}},\mathfrak{S}} P$. From this the proposition follows provided $G$ is connected.   
  
In the case $p \geq 3$ the same argument works if we replace
$(P',Q',F',\dot{F}')$ by the Dieudonn\'e display of $G$ over
$\hat{W}(R)$. But in this case the isomorphism
$\mathbb{D}(G)_{\hat{W}(R)} \cong  P'$ follows from \cite{L2}, Theorem 6.9.  
\endproof

\medskip
All results in this paper are actually independent of the nonconnected
part of the last Proposition. This is explained in the proof of
Corollary \ref{C3}.

In what follows we will not need to keep track of the $\sigma_{\kappa(\ideal{p})^{\perf}}$-linear maps $\sigma_{\kappa(\ideal{p})^{\perf}} \otimes F$ and
thus we will simply call $W(\kappa (\ideal{p})^{\perf})
\otimes_{\varkappa_{\ideal{p}},\mathfrak{S}} P$ the fibre of the Breuil window
$(P,Q,F,\dot{F})$ over $\ideal{p}$.

We often write a Breuil window in the form $(Q,\phi)$ originally proposed by Breuil,
where $\phi$ is the composite of the inclusion $Q\subset P$ with the inverse of the
$\mathfrak{S}$-linear isomorphism $\dot{F}^{\sharp}:Q^{(\sigma)} \cong P$. In this
notation $P$, $F$, and $\dot{F}$ are omitted as they are determined naturally by
$\dot{F}^{\sharp}$ and thus by $\phi$ (see
\cite{VZ}, Section 2). A Breuil window in
this form is characterized as follows: $Q$ is a
finitely generated  
free $\mathfrak{S}$-module and $\phi: Q \rightarrow Q^{(\sigma)}$ is 
a $\mathfrak{S}$-linear map whose cokernel is annihilated
by $E$.  We note that this
implies easily that there exists a unique $\mathfrak{S}$-linear map $\psi: Q^{(\sigma)}
\rightarrow Q$ such that we have
\begin{displaymath}
   \phi \circ \psi = E \id_{Q^{(\sigma)}}, \quad \psi \circ \phi = E
   \id_{Q}. 
\end{displaymath}  
Clearly the datum $(Q,\psi)$ is equivalent to the datum
$(Q,\phi)$. In the notation $(Q,\phi)$, its fibre over $\ideal{p}$ is $W(\kappa
(\ideal{p})^{\perf}) \otimes_{\varkappa_{\ideal{p}},\mathfrak{S}}
Q^{(\sigma)}=W(\kappa (\ideal{p})^{\perf})
\otimes_{\sigma\varkappa_{\ideal{p}},\mathfrak{S}} Q$.

The dual of a Breuil window is defined as follows. Let $M$ be a
$\mathfrak{S}$-module. We set $\hat{M} := \Hom_{\mathfrak{S}} (M,
\mathfrak{S})$. A $\mathfrak{S}$-linear map $M \rightarrow \mathfrak{S}$ defines
a homomorphism $M^{(\sigma)} \rightarrow \mathfrak{S}^{(\sigma)} =
\mathfrak{S}$. This defines a $\mathfrak{S}$--linear map:
\begin{displaymath}
   \hat{M}^{(\sigma)}=\mathfrak{S} \otimes_{\sigma, \mathfrak{S}} 
\Hom_{\mathfrak{S}} (M,\mathfrak{S})
\rightarrow \Hom_{\mathfrak{S}} (M^{(\sigma)},\mathfrak{S})=\widehat{M^{(\sigma)}}.
\end{displaymath}
It is clearly an isomorphism if $M$ is a free $\mathfrak{S}$-module of finite rank 
and therefore also if $M$ is a finitely generated $\mathfrak{S}$-module by a formal 
argument which uses the flatness of $\sigma:\mathfrak{S}\to\mathfrak{S}$.

If $(Q,\phi)$ is a Breuil window we obtain a $\mathfrak{S}$-linear map
\begin{displaymath}
   \hat{\phi}: \hat{Q}^{(\sigma)} = \widehat{(Q^{(\sigma)})}
   \rightarrow \hat{Q}.
\end{displaymath}
More symmetrically we can say that if $(Q,\phi, \psi)$ is a 
Breuil window then $(\hat{Q}, \hat{\psi}, \hat{\phi})$ is a
Breuil window. We call $(\hat{Q}, \hat{\psi}, \hat{\phi})$ the {\it dual Breuil window} of $(Q,\phi, \psi)$. 

Taking the fibre of a Breuil window $(Q, \phi)$ over  $\ideal{p}$ is
compatible with duals as we have:
\begin{equation}\label{Val2e}
   W(\kappa (\ideal{p})^{\perf}) \otimes_{\varkappa_{\ideal{p}},\mathfrak{S}}
   \hat{Q}^{(\sigma)} \cong \Hom_{W(\kappa (\ideal{p})^{\perf})} 
( W(\kappa (\ideal{p})^{\perf}) \otimes_{\varkappa_{\ideal{p}},\mathfrak{S}}
Q^{(\sigma)}, 
W(\kappa (\ideal{p})^{\perf})).
\end{equation}

We recall from \cite{VZ}, Definition 2 and \cite{L1} that a {\it Breuil module}
$(M,\varphi)$ is a pair, where $M$ is a
finitely generated $\mathfrak{S}$-module which is of projective dimension at most
$1$  and which is
annihilated by a power of $p$ and where $\varphi : M \rightarrow
M^{(\sigma)}$ is a $\mathfrak{S}$-linear map whose cokernel is annihilated by $E$.
We note that the map $\varphi$ is always injective (the argument for this is
the same as in \cite{VZ}, Proposition 2 (i)). It follows
formally that there exists a unique $\mathfrak{S}$-linear  map $\vartheta :
M^{(\sigma)}
\rightarrow
M$ such that we have
\begin{displaymath}
   \varphi \circ \vartheta = E \id_{M^{(\sigma)}}, \quad \vartheta \circ \varphi = E
   \id_{M}. 
\end{displaymath}  
We define the {\it dual Breuil module} $(M^{\ast}, \vartheta^{\ast},
\varphi^{\ast})$ of $(M,\varphi,\vartheta)$ by applying the functor $M^{\ast} =
\Ext^1_{\mathfrak{S}}(M, \mathfrak{S})$ in the same manner 
we did for windows. It is easy to see that the $\mathfrak{S}$-module $M^{\ast}$ has
projective dimension at most $1$. The fibre of $(M,\varphi)$ (or of
$(M,\varphi,\vartheta)$) over $\ideal{p}$ is
\begin{equation}\label{Val4e}
   W(\kappa(\ideal{p})^{\perf}) \otimes_{\varkappa_{\ideal{p}}, \mathfrak{S}}
M^{(\sigma)} = 
 W(\kappa(\ideal{p})^{\perf}) \otimes_{\sigma\varkappa_{\ideal{p}}, \mathfrak{S}} M.
\end{equation} 
The duals are again compatible with
taking fibres as we have:
\begin{equation}\label{Val3e}
\begin{array}{l}
W(\kappa (\ideal{p})^{\perf}) \otimes_{\sigma\varkappa_{\ideal{p}},\mathfrak{S}}
M^{\ast}  \cong 
\Ext^1_{W(\kappa (\ideal{p})^{\perf})} (W(\kappa (\ideal{p})^{\perf})
\otimes_{\sigma\varkappa_{\ideal{p}},\mathfrak{S}} M, W(\kappa (\ideal{p})^{\perf}))\\
= \Hom_{W(\kappa (\ideal{p})^{\perf})}(W(\kappa (\ideal{p})^{\perf})
\otimes_{\sigma\varkappa_{\ideal{p}},\mathfrak{S}} M, W(\kappa (\ideal{p})^{\perf})
\otimes_{\mathbb{Z}} \mathbb{Q}). 
\end{array}
\end{equation}

Assume that $p$ annihilates $M$, i.e. $M$ is a module over $\bar{\mathfrak{S}} =
\mathfrak{S}/p\mathfrak{S}$. As $\depth M$ is the same  over either $\mathfrak{S}$
or $\bar{\mathfrak{S}}$ and it is $d$ if $M\neq 0$, we easily get that $M$ is a free
$\bar{\mathfrak{S}}$-module. From this we get that $M^{\ast}=
\Hom_{\bar{\mathfrak{S}}} (M, \bar{\mathfrak{S}})$ (to be compared with the last
isomorphism of  (\ref{Val3e})). Thus in this case the duality works exactly as for
windows.

If $p\geq 3$, it follows from Theorem  \ref{T3} that the category of finite flat
group schemes  over $\Spec R$ is equivalent to the category of Breuil modules (the
argument for this is the same as for \cite{VZ}, Theorem 2). We have a variant of
this for
$p=2$ (cf. \cite{L1}, Theorem 10.7): the category of connected finite flat group schemes  over $\Spec R$ is
equivalent to the category of nilpotent Breuil modules (i.e., of Breuil modules
$(M,\varphi)$ that have the property that the reduction of $\varphi$ modulo the maximal ideal of
$\mathfrak{S}$ is nilpotent in the natural way).

We recall from \cite{VZ}, Subsection 6.1 that the Breuil module of a finite flat
group scheme $H$ over $\Spec R$ is
obtained as follows. By a theorem of Raynaud (see \cite{BBM}, Theorem 3.1.1) we can
represent $H$ as the kernel 
\begin{displaymath}
   0 \rightarrow H \rightarrow G_1 \rightarrow G_2 \rightarrow 0
\end{displaymath}  
of an isogeny $G_1\rightarrow G_2$ of $p$-divisible groups over $\Spec R$. If $p=2$, then we assume that $H$ and $G_1$ are connected. Let $(Q_1,\phi_1)$ and $(Q_2, \phi_2)$ be the
Breuil windows of $G_1$ and $G_2$ (respectively). Then the Breuil module
$(M,\varphi)$ of $H$ is the cokernel
of the induced map $(Q_1, \phi_1) \rightarrow (Q_2, \phi_2)$ in a natural sense. From Proposition \ref{P1} we get that the classical covariant Dieudonn\'e module of $H_{\kappa(\ideal{p})^{\perf}}$ is canonically given by (\ref{Val4e}).

\section{Morphisms between Breuil modules}\label{S1}

In this section, let  $R$ be a complete regular local ring of mixed
characteristic $(0,p)$  with
maximal ideal $\ideal{m}$ and perfect residue class field $k$. We
write  $R = \mathfrak{S}/E\mathfrak{S}$, where $d=\dim R\ge 1$,
$\mathfrak{S}=W(k)[[T_1,\ldots,T_d]]$, and $E = p - h\in \mathfrak{S}$ are as in
Section \ref{S0}. We use the
standard frame $\mathcal L$ of the Definition \ref{D5}. 

Let $e\in \mathbb{N}^*$ be
such that $p \in \ideal{m}^e\setminus \ideal{m}^{e+1}$. It is the
absolute ramification index of $R$. Let $\bar{\ideal{r}}:=(T_1, \ldots, T_d)\subset
\bar{\mathfrak{S}}:=\mathfrak{S}/p\mathfrak{S}=k[[T_1, \ldots, T_d]]$. Let $\bar{h}
\in
\bar{\ideal{r}}\subset \bar{\mathfrak{S}}$ be the reduction modulo $p$ of $h$. The
surjective function $\ord:\bar{\mathfrak{S}}\twoheadrightarrow
\mathbb{N}\cup\{\infty\}$ is such that $\ord
(\bar{\ideal{r}}^i\setminus \bar{\ideal{r}}^{i+1})=i$ for all $i\ge 0$ and
$\ord(0)=\infty$. 

Let $\gr R:=\gr_{\ideal{m}} R$. The obvious isomorphism of graded rings
\begin{equation}\label{Val3e1}
   k[T_1, \ldots, T_d] \rightarrow \gr R
\end{equation}
maps the initial form of $\bar{h}$ to the initial form of $p$. Thus
$e$ is the order of the power series $\bar{h}$.

\begin{lemma}\label{Le1}
We assume that $R$ is such that $p \notin \ideal{m}^p$ (i.e., $e\le p-1$). 
Let $C$ be a $\mathfrak{S}$-module which is annihilated by a power of
$p$. Let $\varphi: C \rightarrow C^{(\sigma)}$ be a
$\mathfrak{S}$-linear map whose cokernel is annihilated by $E$. We
assume that there exists a power series $f \in \mathfrak{S}\setminus
p\mathfrak{S}$ which annihilates $C$. Then we have: 

\medskip
{\bf (a)} If $p \notin \ideal{m}^{p-1}$ (i.e., if $e\le p-2$), then $C=0$.

\smallskip

{\bf (b)} If $e=p-1$, then either $C=0$ or the initial form of $p$  in
$\gr R$ generates an ideal which is a $(p-1)$-th power.


\end{lemma}
{\bf Proof.} It suffices to show that $C=0$ provided either $e\le p-2$
or $e=p-1$ and the initial form of $p$  in $\gr R$ generates an ideal which is not a
$(p-1)$-th power. By the lemma of
Nakayama it suffices to show that $C/pC = 0$. 
It is clear that $\varphi$ induces a $\mathfrak{S}$-linear map $C/pC \rightarrow
(C/pC)^{(\sigma)}$. Therefore we can assume that $C$ is annihilated by
$p$.

Let $u$ be the smallest non-negative integer with the following
property: for each $c \in C$ there exists a power series $g_c \in
\bar{\mathfrak{S}}$ such that $\ord(g_c)\leq u$ and $g_c$ annihilates $c$.

From the existence of $f$ in the annihilator of $C$ we deduce that the 
number $u$ exists. If $C \neq 0$, then we have $u >0$. We will show that the
assumption that $u>0$ leads to a contradiction and therefore we have $C=0$. 

By the minimality of $u$ there exists an element $x \in C$ such that for 
each power series $a$ in the annihilator $\ideal{a} \subset
\bar{\mathfrak{S}}$ of $x$ we have $\ord(a)\geq u$. Consider the
$\bar{\mathfrak{S}}$-linear injection 
\begin{equation}\label{Val1e}
   \bar{\mathfrak{S}}/\ideal{a} \hookrightarrow C
\end{equation}
 which maps $1$ to $x$. Let $\ideal{a}^{(p)} \subset \bar{\mathfrak{S}}$
be the ideal generated by the $p$-th powers of elements in
$\ideal{a}$. Each power series in $\ideal{a}^{(p)}$ has order
$\geq pu$. 

If we tensorize the injection (\ref{Val1e}) by $\sigma:\bar{\mathfrak{S}}\rightarrow
\bar{\mathfrak{S}}$ we obtain a $\bar{\mathfrak{S}}$-linear injection
\begin{displaymath}
   \bar{\mathfrak{S}}/\ideal{a}^{(p)} \cong \bar{\mathfrak{S}}
   \otimes_{\sigma, \bar{\mathfrak{S}}} \bar{\mathfrak{S}}/ \ideal{a}
\hookrightarrow C^{(\sigma)}.
\end{displaymath}
Thus each power series in the annihilator of $1 \otimes x
\in C^{(\sigma)}$ has at least order $pu$. 

On the other hand the cokernel of $\varphi$ is by assumption annihilated
by $\bar{h}$. Thus $\bar{h}(1\otimes x)$ is in the image of
$\varphi$. By the definition of $u$ we find a power series $g \in
\bar{\mathfrak{S}}$ with $\ord(g)\leq u$ which annihilates
$\bar{h}(1\otimes x)$. Thus $g\bar{h} \in \ideal{a}^{(p)}$. We get

\begin{displaymath}
   u + \ord(\bar{h})\geq \ord(g \bar{h}) \geq pu.
\end{displaymath} 
Therefore $e=\ord(\bar{h})\geq (p-1)u\geq p-1$. In the case (a) we obtain a
contradiction which shows that $C = 0$.

In the case $e = p-1$ we obtain a contradiction if $u >1$. 
Assume that $u=1$. As $g\bar{h}\in \ideal{a}^{(p)}$ and as $\ord(g\bar h)=p$, 
there exists a power series $\bar{f} \in \bar{\mathfrak{S}}$ of order $1$ and a non-zero
element $\xi \in k$ such that 
\begin{displaymath}
   g\bar{h} \equiv \xi \bar{f}^p \; \text{mod} \; \bar{\ideal{r}}^{p+1}. 
\end{displaymath}
This shows that the initial forms of $g$ and $\bar{f}$ in the
graded ring $\gr_{\bar{\ideal{r}}} \bar{\mathfrak{S}}$ differ by a constant
in $k$. If we divide the last congruence by $\bar{f}$ we obtain
\begin{displaymath}
     \bar{h} \equiv \xi \bar{f}^{p-1} \; \text{mod} \; \bar{\ideal{r}}^{p}. 
\end{displaymath} 
By the isomorphism (\ref{Val3e1}) this implies that initial form of
$p$ in $\gr R$ generates an ideal which is a $(p-1)$-th
power. Contradiction.\endproof 

\begin{lemma}\label{Le2}
We assume that $d=2$ and we consider the ring $\bar{\mathfrak{S}}=k[[T_1,T_2]] $. Let $\bar{h} \in \bar{\mathfrak{S}}$ be a power series of order $e\in\mathbb{N}^*$. Then the following two statements are equivalent:

\medskip
{\bf (a)} The power series $\bar{h}$ does not belong to the ideal
$\bar{\ideal{r}}^{(p)}+\bar{\ideal{r}}^{2(p-1)}=(T_1^p,
T_2^p,\break T_1^{p-1} T_2^{p-1})$ of $\bar{\mathfrak{S}}$.  

\smallskip
{\bf (b)} If $C$ is a $\bar{\mathfrak{S}}$-module of finite length equipped with a $\bar{\mathfrak{S}}$-linear map $\varphi:C\rightarrow C^{(\sigma)}$ whose
cokernel is annihilated by $\bar h$, then $C$ is zero. 
\end{lemma}

{\bf Proof.} We consider the $\bar{\mathfrak{S}}$-linear map 
$$\tau:k=k[[T_1,T_2]]/(T_1,T_2) \rightarrow k^{(\sigma)}=k[[T_1,T_2]]/(T_1^p,T_2^p)$$ 
that maps $1$ to $T_1^{p-1} T_2^{p-1}$ modulo $(T_1^p,T_2^p)$. If $\bar{h}\in\bar{\ideal{r}}^{(p)}+\bar{\ideal{r}}^{2(p-1)}$, then the cokernel of $\tau$ is annihilated by $\bar h$. From this we get that (b) implies (a).

Before proving the other implication, we first make some general remarks.
The Frobenius endomorphism $\sigma: \bar{\mathfrak{S}} \rightarrow
\bar{\mathfrak{S}}$ is faithfully flat. If $\ideal{\vartriangle}$ is an ideal of $\bar{\mathfrak{S}}$, then with the same notations as
before we have  $\sigma(\ideal{\vartriangle})\bar{\mathfrak{S}} =
\ideal{\vartriangle}^{(p)}$.

Let $C$ be a $\bar{\mathfrak{S}}$-module of finite type. Let
$\ideal{a} \subset \bar{\mathfrak{S}}$ be the annihilator of $C$. Then 
$\ideal{a}^{(p)}$ is the annihilator of $C^{(\sigma)} = \bar{\mathfrak{S}}
\otimes_{\sigma, \bar{\mathfrak{S}}} C$. This is clear for a monogenic
module $C$. If we have more generators for $C$, then we can use the formula
\begin{displaymath}
   \ideal{a}^{(p)} \cap \ideal{b}^{(p)} = (\ideal{a} \cap
   \ideal{b})^{(p)}
\end{displaymath}
which holds for flat ring extensions in general. 

We say that a power series $f \in \bar{\mathfrak{S}}$ of order $u$ is
normalized with respect to $(T_1,T_2)$ if it contains $T_1^u$ with a
non-zero coefficient. This definition makes sense with respect to any
regular system of parameters $\tilde T_1,\tilde T_2$ of
$\bar{\mathfrak{S}}$.

We now assume that (a) holds and we show that (b) holds.  Thus the $\bar{\mathfrak{S}}$-module $C$ has finite length and the
cokernel of $\varphi:C\to C^{(\sigma)}$ is annihilated by $\bar h$.  Let $k^\prime$ be an infinite perfect field that contains $k$. Let $\bar{\mathfrak{S}}^\prime:=k^\prime[[T_1,T_2]]$. To show that $C=0$, it suffices to show that $\bar{\mathfrak{S}}^\prime\otimes_{\bar{\mathfrak{S}}} C=0$. Thus by replacing the role of $k$ by the one of $k^\prime$, we can assume that $k$ is infinite. This assumption implies that for almost all $\lambda \in k$, the
power series $\bar h$ is normalized with respect to $(T_1, T_2 + \lambda
T_1)$. By changing the regular system of parameters $(T_1,T_2)$ in
$\bar{\mathfrak{S}}$, we can assume that $\bar h$ is normalized with
respect to $(T_1,T_2)$. By the Weierstra{\ss} preparation theorem we
can assume that $\bar h$ is a Weierstra{\ss} polynomial (\cite{Bou}, Chapter 7,
Section 3, number 8). Thus we can write

\begin{displaymath}
   \bar{h} = T_1^e + a_{e-1}(T_2)T_1^{e-1} + \ldots +a_1(T_2)T_1 + a_0(T_2),
\end{displaymath}
where $a_0(T_2),\ldots,a_{e-1}(T_2) \in T_2k[[T_2]]$. 

We note that the assumption (a) implies that $e \leq 2p -3$.   
 
Let $u$ be the minimal non-negative integer such that there exists a power series of the form $T_1^u + g \in \bar{\mathfrak{S}}$, with $g \in
T_2\bar{\mathfrak{S}}$, for which we have $(T_1^u + g) C = 0$. By our
assumptions, such a non-negative integer $u$ exists.  

We will show that the assumption that $C \neq 0$ leads to a contradiction. This assumption implies that $u\ge 1$. The annihilator $\ideal{a}$ of the
module $C$ has a set of generators of the following form:

\begin{displaymath}
   T_1^{u_i} + g_i, (\text{with} \; i = 1, \ldots, l), \quad g_i \;
   (\text{with} \; i = l+1, \ldots, m),
\end{displaymath}
where $u_i \geq u$ for $i = 1, \ldots, l$ and $g_i \in
T_2\bar{\mathfrak{S}}$ for $i=1, \ldots, m$.  

As the cokernel of $\varphi$ is annihilated by $\bar{h}$ we find
that $(T_1^u + g)\bar{h}C^{(\sigma)} = 0$.
As the annihilator of $C^{(\sigma)}$ is generated by the elements 

\begin{equation}\label{Val5e}
   T_1^{pu_i} + g_i^p\; (\text{with} \; i = 1, \ldots, l), \quad g_i^p \;
(\text{with} \; i
   = l+1, \ldots, m),
\end{equation}
we obtain the congruence
\begin{equation}\label{Val5e1}
  (T_1^u + g)\bar{h} \equiv 0 \; \mod (T_1^{up}, T_2^p). 
\end{equation}
We consider this congruence modulo $T_2$. We have $g \equiv 0 \mod
(T_2)$ and $\bar h \equiv T_1^e\;  \mod (T_2)$ because $\bar h$ is a Weierstra{\ss}
polynomial. This proves that
\begin{displaymath}
   T_1^uT_1^e \equiv 0 \mod (T_1^{up}).
 \end{displaymath}
But this implies that $u + e \geq up$. If $u\geq 2$, then $e\geq up-u\geq 2p-2$ and this contradicts the inequality $e \leq 2p -
3$. Therefore we can assume that $u = 1$. 

By replacing $(T_1,T_2)$ with $(T_1+g,T_2)$, without loss of generality we can assume that $T_1 C = 0$ and (cf. (\ref{Val5e1}) and the equality $u=1$) that 
\begin{displaymath}
   T_1 \bar{h} \equiv 0 \; \mod \bar{\ideal{r}}^{(p)}.
\end{displaymath}
This implies that up to a unit in $\bar{\mathfrak{S}}$ we can assume
that $\bar{h}$ is of the form:

\begin{displaymath}
   \bar{h} \equiv T_2^sT_1^{p-1} + \sum_{i=0}^{\infty}
   T_2^{i+p}\delta_i(T_1) \; \mod (T_1^p),
\end{displaymath}
where $0\leq s \leq p-2$ and where each $\delta_i\in k[T_1]$ has degree at most $p-2$. 

Let now $v$ be the smallest natural number such that $T_2^v C = 0$. 
Then the annihilator $\ideal{a}$ of $C$ is generated by $T_1, T_2^v$ and the
annihilator $\ideal{a}^{(p)}$ of $C^{(\sigma)}$ is generated by $T_1^p, T_2^{pv}$. 
As the cokernel of $\varphi$ is annihilated by $\bar h$, we find that $T_2^v\bar h C^{(\sigma)}=0$. Thus we obtain the congruence
\begin{displaymath}
   T_2^v(T_2^sT_1^{p-1} + \sum_{i=0}^{\infty} T_2^{i+p}\delta_i(T_1))
   \equiv 0 \; \mod (T_1^p, T_2^{pv}).
\end{displaymath}
But this implies $v + s \geq pv$. Thus $s\ge (p-1)v\ge p-1$ and (as $0\leq s \leq p-2$) we reached a contradiction. Therefore $C=0$ and thus (a) implies (b).\endproof

\begin{prop}\label{P2}
We assume that $p \notin \ideal{m}^p$ (i.e., $e\leq p-1$). If $p\in\ideal{m}^{p-1}$,
then
we also assume that the ideal generated by the initial form of $p$ in
$\gr R$ is not a $(p-1)$-th power (thus $p>2$). 

We consider a morphism $\alpha: (M_1, \varphi_1) \rightarrow (M_2, \varphi_2)$ of
Breuil modules for the standard frame $\mathcal{L}$. 
Let $ \ideal{p}$ be a prime ideal
of $R$ which contains $p$. We consider the
$W(\kappa(\ideal{p})^{\perf})$-linear map obtained from $\alpha$ by
base change 
\begin{equation}\label{Val6e}
W(\kappa(\ideal{p})^{\perf}) \otimes_{\sigma \varkappa_{\ideal{p}}, \mathfrak{S}}
M_1 
\rightarrow  
W(\kappa(\ideal{p})^{\perf}) \otimes_{\sigma \varkappa_{\ideal{p}}, \mathfrak{S}}
M_2.
\end{equation}
Then the following two properties hold:

\medskip
{\bf (a)} If (\ref{Val6e}) is surjective, then the $\mathfrak{S}$-linear map $M_1
\rightarrow M_2$ is surjective.

{\bf (b)} If (\ref{Val6e}) is
injective, then the $\mathfrak{S}$-linear map $M_1 \rightarrow M_2$ is injective and
its cokernel
is a $\mathfrak{S}$-module of projective dimension at most $1$.
\end{prop}
{\bf Proof.}  We prove (a). We denote by $\tilde{\ideal{p}}$ the ideal of
$\mathfrak{S}$ which corresponds to $\ideal{p}$ via the isomorphism
$\mathfrak{S}/(E,p) \cong R/(p)$. It follows from the lemma of 
Nakayama that $(M_1)_{\tilde{\ideal{p}}} \rightarrow
(M_2)_{\tilde{\ideal{p}}}$ is a surjection. Let us denote by $(C,
\varphi)$ the cokernel of $\alpha$. We conclude that $C$ is annihilated 
by an element $f \notin \tilde{\ideal{p}} \supset p\mathfrak{S}$.
Therefore we conclude from Lemma \ref{Le1} that $C = 0$. Thus (a) holds.

Part (b) follows from (a) by duality
 as it is compatible with taking fibres (\ref{Val3e}). Indeed for 
the last assertion it is enough to note that the kernel of the
surjection $M_2^{\ast} \twoheadrightarrow M_1^{\ast}$ has clearly projective
dimension at most $1$. \endproof

\medskip\medskip
\noindent
{\bf Remark.} If $e=p-1$ and the ideal generated by the initial form of $p$ in
$\gr R$ is a $(p-1)$-th power, then the Proposition \ref{P2} is not true in general.
This is so as there  exist non-trivial homomorphisms
$(\mathbb{Z}/p\mathbb{Z})_R\rightarrow
\pmb{\mu}_{p,R}$ for suitable such $R$'s.

\begin{prop}\label{P3}
We assume that $R$ has dimension $d = 2$ and that $h \in
\mathfrak{S}$ does not belong to the ideal
$(p,T_1^p, T_2^p, T_1^{p-1}T_2^{p-1})$ of
$\mathfrak{S}$. Let $\alpha: (M_1, \varphi_1) 
\rightarrow (M_2, \varphi_2)$ be a morphism of Breuil modules for the standard frame
$\mathcal{L}$. We assume that for each prime ideal $\ideal{p}$ of $R$ with $p \in
\ideal{p} \neq \ideal{m}$, the $W(\kappa(\ideal{p})^{\perf})$-linear map obtained
from $\alpha$ by
base change
\begin{equation}\label{Val7e}
W(\kappa(\ideal{p})^{\perf}) \otimes_{\sigma \varkappa_{\ideal{p}}, \mathfrak{S}}
M_1 
\rightarrow  
W(\kappa(\ideal{p})^{\perf}) \otimes_{\sigma \varkappa_{\ideal{p}}, \mathfrak{S}}
M_2 
\end{equation}
is surjective (resp. is injective).

Then the $\mathfrak{S}$-linear map $M_1 \rightarrow M_2$ is surjective (resp. is
injective and its cokernel
is a $\mathfrak{S}$-module of projective dimension at most $1$).
\end{prop}
{\bf Proof.}   As in the last proof we only have to treat the case where 
the maps (\ref{Val7e}) are surjective. We consider the cokernel
$(C,\varphi)$ of $\alpha$. As in the last proof we will argue that $C=0$ but with the role of Lemma
\ref{Le1} being replaced by Lemma \ref{Le2}. It suffices to show that $\bar C:=C/pC$ is zero. Let $\bar\varphi:\bar C\to \bar C^{(\sigma)}$ be the map induced naturally by $\varphi$.

Lemma \ref{Le2} is applicable if we verify that $\bar C$ is a
$\bar{\mathfrak{S}}$-module of finite length. 

We denote by 
$\bar{\ideal{p}}$ the ideal of $\bar{\mathfrak{S}} =
\mathfrak{S}/p\mathfrak{S}$ which corresponds to $\ideal{p}$ via the 
isomorphism $\bar{\mathfrak{S}}/(\bar{h}) \cong R/pR$. It follows from
(\ref{Val7e}) by the lemma of Nakayama that for each  prime ideal
$\bar{\ideal{p}} \supset \bar{h} \bar{\mathfrak{S}}$ different from
the maximal ideal of $\bar{\mathfrak{S}}$ the maps 
$\alpha_{\ideal{p}}:(M_1)_{\bar{\ideal{p}}} \rightarrow (M_2)_{\bar{\ideal{p}}}$ are
surjective. Using this one constructs inductively a regular sequence
$f_1, \ldots, f_{d-1}, \bar{h}$ in the ring $\bar{\mathfrak{S}}$ such 
that the elements $f_1, \ldots, f_{d-1}$ annihilate $\bar{C}$. As $\bar{C}$ is
a finitely generated module over the $1$-dimensional local ring
$\bar{\mathfrak{S}}/(f_1,
\ldots, f_{d-1})$, we get that $\bar{C}[1/\bar{h}]$ is a module of finite
length over the regular ring $A = \bar{\mathfrak{S}}[1/\bar{h}]$. 
If we can show that $\bar{C}[1/\bar{h}] = 0$, then it
follows that $\bar{C}$ is of finite length.  

As we are in characteristic
$p$, the Frobenius $\sigma$ acts on the principal ideal domain $A =
\bar{\mathfrak{S}}[1/\bar{h}]$. By the definition of a Breuil module
the maps $\varphi_i [1/\bar{h}]$ for $i = 1,2$ become
isomorphisms. Therefore $\varphi$ gives birth to an isomorphism: 
\begin{equation}\label{Val8e}
   \bar{\varphi}[1/\bar{h}]:\bar{C}[1/\bar{h}] \rightarrow
(\bar{C}[1/\bar{h}])^{(\sigma)}. 
\end{equation}
As $A$ is regular of dimension $d-1$, for each $A$-module $\ddag$ of finite length we
have 
\begin{displaymath}
   \length \ddag^{(\sigma)} = p^{d-1} \length\ddag.
\end{displaymath}
We see that the isomorphism (\ref{Val8e}) is only possible if 
$\bar{C}[1/\bar{h}] = 0$. Thus $\bar{C}$ has finite length and therefore
from Lemma \ref{Le2} we get that $\bar C=0$. This implies that $C=0$.\endproof

\section{Extending epimorphisms and monomorphisms}

In this section let $R$ be a regular local ring of mixed characteristic
$(0,p)$ with maximal ideal $\ideal{m}$ and residue class field
$k$. Let $K$ be the field of fractions of $R$.

\subsection{Complements on Raynaud's work}
We first reprove Raynaud's result \cite{R2}, Corollary 3.3.6 by the methods of the
previous sections. We state it in a slightly different form.

\begin{prop}\label{P4}
We assume that $p \notin \ideal{m}^{p-1}$ (thus $p>2$). Let $H_1$ and $H_2$ be
finite flat group schemes over $\Spec R$. 
Let  $ \beta: H_1 \rightarrow H_2$ be a homomorphism, which induces
an epimorphism (resp. monomorphism) $H_{1,K} \rightarrow H_{2,K}$ between generic 
fibres.  Then $\beta$ is an epimorphism (resp. monomorphism).  
\end{prop}
{\bf Proof.} We prove only the statement about epimorphisms because
the case of a monomorphism follows by Cartier duality.  
It is enough to show that the homomorphism $H_1 \rightarrow
H_2$ is flat. By the fibre criterion of flatness it is enough to show
that the homomorphism $\beta_k:H_{1,k}
\rightarrow H_{2,k}$ between special fibres is an epimorphism. To see this we can assume that
$R$ is a complete local ring with algebraically closed residue class field
$k$. 

We write $R \cong \mathfrak{S}/(p - h)$, with $\mathfrak{S}=W(k)[[T_1,\ldots,T_d]]$.
Then the reduction $\bar{h} \in k[[T_1,\ldots,T_d]]$ of $h$ modulo $p$ is a power
series of order $e < p-1$. By
Noether normalization theorem we can assume that $\bar{h}$ contains the monom
$T^{e}_1$. By replacing $R$ with $R/(T_2, \ldots, T_d)$, we can assume that
$R$ is one-dimensional. 

We consider the morphism  
\begin{displaymath}
   (M_1,\varphi_1) \rightarrow (M_2, \varphi_2)
\end{displaymath}
of Breuil modules associated to $\beta$.
Let $(C,\varphi)$ be its cokernel. Let $\bar{C}:=C/C_0$, where $C_0$ is the
$\mathfrak{S}$-submodule of $C$ whose elements are annihilated by a power
of the maximal ideal $\ideal{r}$ of $\mathfrak{S}$. The $\mathfrak{S}$-linear map
$\varphi $ factors as
\begin{displaymath}
   \bar{\varphi}:\bar{C} \rightarrow \mathfrak{S} \otimes_{\sigma,\mathfrak{S}}\bar{C}
\end{displaymath} 
and the cokernel of $\bar{\varphi}$ is annihilated by $\bar h$. The maximal ideal $\ideal{r}$ of $\mathfrak{S}$ is not associated to
$\mathfrak{S} \otimes_{\sigma,\mathfrak{S}}\bar{C}$. Thus either $\bar{C}=0$ or
$\depth \bar{C}\ge 1$. Therefore the $\mathfrak{S}$-module $\bar{C}$ is of
projective dimension at most $1$. As $\bar{C}$ is annihilated by a power of $p$, we
conclude that
$(\bar{C},\bar{\varphi})$
is the Breuil module of a finite flat group scheme $D$ over $\Spec R$. We
have induced homomorphisms
\begin{displaymath}
   H_1 \rightarrow  H_2 \rightarrow D. 
\end{displaymath}
The composition of them is zero and the second homomorphism is an
epimorphism because it is so after base change to $k$. As $H_{1,K}
\twoheadrightarrow H_{2,K}$ is an epimorphism we conclude that $D_K = 0$.
But then $D  = 0$ and the Breuil module $(\bar{C},\bar{\varphi})$ is zero as
well. We conclude by Lemma \ref{Le1} that $C = 0$.\endproof

\medskip
The next proposition is proved in \cite{R2}, Remark 3.3.5 in the case of biconnected finite flat group schemes $H$ and $D$.

\begin{prop}\label{P11}
Let $R$ be a discrete valuation ring of mixed characteristic
$(0,p)$ and index of ramification $p-1$; we have $K = R[1/p]$. 
Let $\beta: H \rightarrow D$ be a homomorphism of residually connected finite flat
group 
schemes over $\Spec R$ which induces an isomorphism (resp. epimorphism) over $\Spec
K$. 

Then $\beta$ is an isomorphism (resp. epimorphism).
\end{prop}
{\bf Proof.} It is enough to show the statement about
isomorphisms. Indeed, assume that $\beta_K$ is an
epimorphism. Consider the schematic closure $H_1$ of the kernel of $\beta_K$ in $H$.
Then $H/H_1 \rightarrow D$ is an
isomorphism.

Therefore we can assume that $\beta_K$ is an isomorphism. By extending $R$
we can assume that $R$ is complete and that $k$ is
algebraically closed. By considering the Cartier dual homomorphism
$\beta^{\text{t}}:D^{\text{t}}\rightarrow H^{\text{t}}$ one easily reduces the
problem to the case when $H$
is biconnected. As the case when $D$ is also
biconnected is known (cf. \cite{R2}, Remark 3.3.5), one can easily reduce to the case when $D$ is
of multiplicative type. Then $D$ contains $\pmb{\mu}_{p,R}$ as a closed subgroup
scheme. We consider the schematic closure $H_1$ of $\pmb{\mu}_{p,K}$ in $H$. Using
an induction on the order of $H$ it is enough to show that
the natural homomorphism $\beta_1:H_1\rightarrow \pmb{\mu}_{p,R}$ is an isomorphism.
This follows from \cite{R2}, Proposition 3.3.2 $3^\circ$. For the sake of
completeness we reprove this in the spirit of the paper.

We write $R = \mathfrak{S}/(E)$, where $E \in \mathfrak{S}=W(k)[[T]]$ is an
Eisenstein polynomial of degree $e = p-1$. 
The Breuil window of the $p$-divisible group $\pmb{\mu}_{p^{\infty},R}$ is given by
\begin{displaymath}
\begin{array}{lcr}
\mathfrak{S} & \rightarrow & \mathfrak{S}^{(\sigma)}\cong \mathfrak{S}\\
 f & \mapsto & Ef.
\end{array}
\end{displaymath} 
The Breuil module $(N, \tau)$ of $\pmb{\mu}_{p,R}$ is the kernel of the
multiplication by $p:\pmb{\mu}_{p^{\infty},R}\rightarrow  \pmb{\mu}_{p^{\infty},R}$
and therefore it can be identified with 
\begin{displaymath}
  \begin{array}{lcr}
N = k[[T]] & \rightarrow & N^{(\sigma)} = k[[T]]^{(\sigma)}\cong k[[T]]\\
 f & \mapsto & T^ef.
\end{array} 
\end{displaymath} 
Let $(M,\varphi)$ be the Breuil module of $H_1$. As $H_1$ is of height
1 we can identify $M = k[[T]]$. Then $\varphi:M\to M^{(\sigma)}\cong M$ is the
multiplication by a power series $g\in k[[T]]$ of order $\ord(g)\leq e=p-1$. To the
homomorphism $\beta_1:H_1\rightarrow \pmb{\mu}_{p,R}$ corresponds a morphism
$\alpha_1:(M,\varphi)\rightarrow (N,\tau)$ that maps $1\in M$ to some element $a\in
N$. We
get a commutative diagram:
\begin{displaymath}
   \xymatrix{
k[[T]] \ar[r]^{a} \ar[d]_g & k[[T]] \ar[d]^{T^e}\\
k[[T]] \ar[r]^{a^p} & k[[T]].
}
\end{displaymath}
 We obtain the equation $ga = a^pT^{p-1}$. As $a \neq 0$ this is
 only possible if 
\begin{displaymath}
   \ord(g) = (p-1)(\ord(a)+ 1).
\end{displaymath} 
As $\ord(g)\leq p-1$, we get that $\ord(a) = 0$ and $\ord(g) = p-1$. As $\ord(a) =
0$, both $\alpha_1$ and $\beta_1$ are isomorphisms. Thus $\beta:H\rightarrow D$ is an isomorphism. \endproof

\subsection{Basic extension properties}\label{S2.5}
The next three results will be proved for $p=2$ under certain residually connectivity assumptions. But in Subsection \ref{S2} below we will show how
these results hold under no residually connectivity assumptions even if $p=2$.

\begin{prop}\label{P5}
We assume that $\dim R \geq 2$, that $p \notin \ideal{m}^{p}$, and that 
the initial form of $p$ in $\gr R$ generates an ideal which is not a $(p-1)$-th
power. Let $\beta: H_1 \rightarrow H_2$ be a homomorphism of finite flat group
schemes 
over $\Spec R$ which induces an epimorphism (resp. a monomorphism) over
$\Spec K$. If $p=2$, we assume as well that $H_1$ and $H_2$ (resp. that the Cartier
duals
of $H_1$ and $H_2$) are residually connected.

Then $\beta$ is an epimorphism (resp. a monomorphism).
\end{prop}
{\bf Proof.} As before we can assume that $R$ is a complete regular
local ring. We can also restrict our attention to the case of epimorphisms. 

By Proposition \ref{P4} we can assume that $p \in \ideal{m}^{p-1}$. Let
$\ideal{p}$ be a minimal prime ideal which contains $p$. We show that the assumption that $p \in \ideal{p}^{p-1}$ leads to a contradiction. Then we can
write $p = uf^{p-1}$, where $u,f \in R$ and $f$ is a generator of the
prime ideal $\ideal{p}$. It follows from our assumptions that $u \notin \ideal{m}$
and that $f \in \ideal{m} \setminus \ideal{m}^2$.  Therefore the initial form of $p$
in $\gr R$ generates an ideal which is a $(p-1)$-th power. Contradiction.

We first consider the case when $k$ is perfect. Let $\alpha: (M_1, \varphi_1)
\rightarrow (M_2, \varphi_2)$ denote also the morphism of Breuil modules associated
to $\beta:H_1\rightarrow H_2$. As $p \notin
\ideal{p}^{p-1}$, we can apply Proposition \ref{P4} to the ring
$R_{\ideal{p}}$. It follows that $\beta$ induces an epimorphism over the spectrum of
the residue
class field $\kappa(\ideal{p})$ of $\ideal{p}$ and thus also of its perfect hull
$\kappa(\ideal{p})^{\perf}$. From this and the end of Section \ref{S0} we get that the hypotheses of Proposition \ref{P2} hold for $\alpha$. We conclude that $\alpha$ is
an epimorphism and thus $\beta$ is also an epimorphism. 

Let $R \rightarrow R'$ be a faithfully flat extension of noetherian local
rings, such that $\ideal{m}R'$ is the maximal ideal of $R'$ and the extension of
residue class fields $k\hookrightarrow k'$ is radical. We consider the
homomorphism of polynomial rings $\gr R \rightarrow \gr R'$. As $\gr R$ and $\gr R'$
are unique factorization domains, it is easy to see that the
condition that the initial form of $p$ is not a $(p-1)$-th power is
stable by the extension $R \rightarrow R'$. But $\beta$ is an epimorphism if and
only if $\beta_{R^\prime}$ is so. Therefore we can assume
that the residue class field of $R$ is perfect and this case was already proved.
\endproof

\begin{prop}\label{P6}
We assume that $\dim R = 2$. Let $U = \Spec R \setminus \{\ideal{m} \}$. 
We also assume that the following technical condition holds:

\medskip
{\bf (*)} there exists a faithfully flat local $R$-algebra $\hat{R}$
which is complete, which has an algebraically closed
residue class field $k$, and which
has a presentation
$\hat{R} = \mathfrak{S}/(p - h)$ where $\mathfrak{S}=W(k)[[T_1,T_2]]$ and where $h\in (T_1,T_2)$ does not belong to the ideal
$(p,T_1^p, T_2^p,T_1^{p-1}T_2^{p-1})$. 
 
\medskip
Let $\beta: H_1 \rightarrow H_2$ be a homomorphism of finite flat group
schemes  over $\Spec R$. We also assume that for each geometric point $\Spec L
\rightarrow  U$ such that $L$ has characteristic $p$ the homomorphism
$\beta_L$ is an epimorphism (resp. a monomorphism); thus $\beta_U:H_{1,U}\rightarrow
H_{2,U}$ is an epimorphism (resp. a monomorphism). If $p=2$, we assume
that $H_1$ and $H_2$ (resp. that the Cartier duals of $H_1$ and 
$H_2$) are residually connected.

Then $\beta $ is an epimorphism (resp. a
monomorphism).
\end{prop}
{\bf Proof.} We have $\dim \hat{R}/\ideal{m}\hat{R}=0$ and thus
$\Spec(\hat{R})\setminus (\Spec(\hat{R})\times_{\Spec R} U)$ is the closed point of
$\Spec(\hat{R})$. Based on this we can assume that $R=\hat{R}$. Thus the proposition
follows from Proposition \ref{P3} in the same
way Proposition \ref{P5} followed from Propositions \ref{P4} and \ref{P2}.  \endproof

\begin{cor}\label{C2}
We assume that the assumptions on $R$ of either Proposition \ref{P5} or
Proposition \ref{P6} are satisfied (thus $\dim R\geq 2$). We consider a complex
\begin{displaymath}
   0 \rightarrow H_1 \rightarrow H_2 \rightarrow H_3 \rightarrow 0
\end{displaymath}
of finite flat group schemes over $\Spec
R$ whose restriction to $U = \Spec R \setminus \{\ideal{m} \}$ is a short exact sequence. If $p=2$, then we assume that
$H_2$ and $H_3$ are residually connected. Then $0 \rightarrow H_1 \rightarrow H_2
\rightarrow H_3 \rightarrow 0$ is a short exact sequence.
\end{cor}
{\bf Proof.} Propositions \ref{P5} and \ref{P6} imply that $H_2 \rightarrow
H_3$ is an epimorphism. Its kernel is a finite flat group scheme isomorphic to $H_1$
over $U$ and thus (as we have $\dim R\geq 2$) it is isomorphic to $H_1$.\endproof

\subsection{D\'evissage properties}\label{S2}

In this subsection we explain how one can get the results of Subsection \ref{S2.5} under no residually connectivity assumptions for $p=2$. 

We assume that $\dim R=2$ and that $p$ is arbitrary.  Let
$U=\Spec R\setminus\{\ideal{m}\}$. Let $\mathcal{V}$ be a locally free 
$\mathcal{O}_U$-module of finite rank over $U$. The $R$-module $H^0(U,
\mathcal{V})$ is free of finite rank. 
Using this it is easy to see that each finite flat group scheme over
$U$ extends uniquely to a finite flat group scheme over $\Spec
R$. This is clearly an equivalence between the category of finite flat group
schemes over $U$ and the category of finite flat group schemes over $\Spec R$. The
same holds if we restrict to finite flat group schemes annihilated by $p$. 

Next we also assume that $R$ is complete. Each finite flat group scheme $H$ over
$\Spec R$ is canonically an extension
\begin{displaymath}
   0 \rightarrow H^{\circ} \rightarrow H \rightarrow H^{\text{\'et}} \rightarrow 0, 
\end{displaymath}
where $H^{\circ}$ is connected and $H^{\text{\'et}}$ is \'etale over
$\Spec R$. In particular a homomorphism from a connected finite flat group scheme over
$\Spec R$ to an \'etale finite flat group scheme over
$\Spec R$ is zero. From this one easily checks that if $H_1 \rightarrow H_2$ is a
homomorphism of
finite flat group schemes over
$\Spec R$ which is an epimorphism over $U$ and if $H_1$ is connected, then $H_2$ is
connected as well.

\begin{lemma}\label{Le3}
We assume that $\dim R=2$ and that $R$ is complete.  Then the following four
statements are equivalent:
\begin{itemize}
\item[(a)] Each short exact sequence of finite flat group schemes over
  $U$ extends uniquely to a short exact sequence of finite flat group
  schemes over $\Spec R$.  
\item[(b)] Let $H_1$ and $H_2$ be connected finite flat group
  schemes over $\Spec R$. A homomorphism $H_1 \rightarrow H_2$ over $\Spec R$ is an
epimorphism
if its
  restriction to $U$ is an epimorphism. 
\item[(c)] Let $H_1$ and $H_2$ be connected finite flat group
  schemes over $\Spec R$ which are annihilated by $p$. A homomorphism $H_1
  \rightarrow H_2$ over $\Spec R$ is an epimorphism if its 
  restriction to $U$ is an epimorphism.
\item[(d)]  The regular ring $R$ is $p$-quasi-healthy.
\end{itemize}
\end{lemma}
{\bf Proof.} It is clear that (a) implies (b). We show that (b) implies (a). Let 
\begin{displaymath}
   0 \rightarrow H_1 \rightarrow H_2 \rightarrow H_3 \rightarrow 0
\end{displaymath} 
be a complex of finite flat group schemes over $\Spec R$ whose restriction to $U$ is
a short
exact sequence. 
It suffices to show that $\beta: H_1\rightarrow H_2$ is a monomorphism. Indeed, in
this case we
can form the quotient group scheme $H_2/H_1$ and the homomorphism $H_2/H_1
\rightarrow H_3$ is an isomorphism as its restriction to $U$ is so. 

We check that the homomorphism
\begin{displaymath}
   \beta^{\circ} : H_1^{\circ} \rightarrow H_2^{\circ}
\end{displaymath}
is a monomorphism. Let $H_4$ be the finite flat group scheme over $R$
whose restriction to $U$ is $H_{2,U}^{\circ}/H_{1,U}^{\circ}$. We have
a complex $0\rightarrow  H_1^{\circ} \rightarrow H_2^{\circ}\rightarrow
H_4\rightarrow 0$
whose restriction to $U$ is exact. We conclude that $H_4$ is
connected. Thus we have an epimorphism
$H_2^{\circ}\twoheadrightarrow H_4$ (as we are assuming that (b) holds) whose kernel
is $H_1^{\circ}$. Therefore (b) implies that $\beta^{\circ}$ is a monomorphism.  

As $\beta^{\circ}$ is a monomorphism, it suffices to show that the induced
homomorphism $\bar{\beta} : H_1/H_1^{\circ} \rightarrow H_2/H_1^{\circ}$
is a monomorphism. In other words, without loss of generality 
we can assume that $H_1=H_1^{\text{\'et}}$ is \'etale.

Let $H_1''$ be the kernel of $H_2^{\text{\'et}} \rightarrow
H_3^{\text{\'et}}$; it is a finite \'etale group scheme over $\Spec R$. Let $H_1'$
be the
kernel of $H_1
\rightarrow H_1''$.  It suffices to show that $H_1' \rightarrow
H_2^{\circ}$ is a monomorphism. Therefore we can also assume that $H_2$ is
connected. This implies that $H_3$ is connected.  As we are assuming that (b) holds,
$H_2\rightarrow H_3$ is an epimorphism. Its kernel is $H_1$ and therefore $\beta:H_1
\rightarrow H_2$ is a monomorphism. Thus (b) implies (a).

We show that (c) implies (b). The last argument shows that a short exact
sequence of finite flat group schemes over $U$ annihilated by $p$ extends
to a short exact sequence of finite flat group schemes over $\Spec R$. 
We start with a homomorphism $H_2 \rightarrow H_3$ between connected finite flat
group schemes over $\Spec R$ which induces an epimorphism over
$U$. We extend the kernel of  $H_{2,U} \rightarrow H_{3,U}$ to a
finite flat group scheme $H_1$ over $\Spec R$. Then we find a complex
\begin{displaymath}
   0 \rightarrow H_1 \rightarrow H_2 \rightarrow H_3 \rightarrow 0
\end{displaymath} 
whose restriction to $U$ is a short exact sequence. To show that $H_2\rightarrow
H_3$ is an
epimorphism it is equivalent to show that $H_1\rightarrow H_2$ is a monomorphism. As
$H_{1,U}$ has a composition
series whose factors are annihilated by $p$, we easily reduce to the
case where $H_1$ is annihilated by $p$. We embed $H_2$ into a
$p$-divisible group $G$ over $\Spec R$. To check that $H_1\rightarrow H_2$ is a
monomorphism, it suffices to show that $H_1 \rightarrow
G[p]$ is a monomorphism. But this follows from the second sentence of this paragraph. 

It is clear that (a) implies (d). We are left to show that (d) implies
(c). It suffices to show that a homomorphism $\beta: H_1 \rightarrow H_2$ of
finite flat group schemes over $\Spec R$ annihilated by $p$ is a
monomorphism if its restriction $\beta_U:H_{1,U}\rightarrow H_{2,U}$ is a
monomorphism. 
We embed $H_2$ into a $p$-divisible group $G$ over $\Spec R$. The quotient
$G_U/H_{1,U}$ is a $p$-divisible group over $U$ which extends to a
$p$-divisible group $G'$ over $\Spec R$ (as we are assuming that (d) holds). The
isogeny $G_U \rightarrow G'_U$
extends to an isogeny $G\rightarrow G'$. Its kernel is a finite flat
group scheme and therefore isomorphic to $H_1$. We obtain a
monomorphism $H_1 \rightarrow G$. Thus $\beta:H_1
\rightarrow H_2$ is a monomorphism, i.e. (d) implies
(c).   
\endproof

\begin{cor}\label{C3}
Propositions \ref{P5} and \ref{P6} (and
thus implicitly Corollary \ref{C2}) hold without any residually connectivity assumption for $p=2$.
\end{cor}
{\bf Proof.} We can assume
that $R$ is complete. By considering an epimorphism $R\twoheadrightarrow R^\prime$ with
$R^\prime$ regular of dimension $2$, we can also assume that $\dim
R=2$. Based on Lemma \ref{Le3}, it suffices to prove Propositions
\ref{P5} and \ref{P6} in the case when connected finite flat group schemes
are involved. But this case was already proved in Subsection \ref{S2.5}. This also shows that for $p \geq 3$ we can restrict to the
connected case and avoid to apply Proposition \ref{P1} for
nonconnected $p$-divisible groups.\endproof

\begin{cor}\label{C4} 
We assume that $k$ is perfect and that $R=W(k)[[T_1,T_2]]/(p-h)$ with $h\in
(T_1,T_2)$. Let $\bar h\in k[[T_1,T_2]]$ be the reduction of $h$ modulo $p$. Then
the fact that $R$ is or is not $p$-quasi-healthy depends only on the orbit of the
ideal $(\bar h)$ of $k[[T_1,T_2]]$ under automorphisms of $k[[T_1,T_2]]$.
\end{cor}  
{\bf Proof.} The category of Breuil modules associated to connected finite flat
group schemes over $\Spec R$ annihilated by $p$ is equivalent to the category of
pairs $(M,\varphi)$, where $M$ is a free $k[[T_1,T_2]]$-module of finite rank and
where $\varphi:M\rightarrow M^{(\sigma)}$ is a $k[[T_1,T_2]]$-linear map whose
cokernel is
annihilated by $\bar h$ and whose reduction modulo the ideal $(T_1,T_2)$ is
nilpotent in the natural sense. The last category depends only on the orbit of
$(\bar h)$ under automorphisms of $k[[T_1,T_2]]$. The corollary follows from the
last two sentences and the equivalence of (c) and (d) in Lemma \ref{Le3}.\endproof

\subsection{The  $p$-quasi-healthy part of Theorem \ref{T1}}\label{S3} In this
subsection we show that if $R$ is as in Theorem \ref{T1} for $d=2$, then $R$ is
$p$-quasi-healthy. It follows from the definition of a $p$-divisible group
and the uniqueness part of the first paragraph of Subsection \ref{S2}, that it is
enough to show that a complex $0\rightarrow H_1\rightarrow H_2\rightarrow
H_3\rightarrow 0$ of finite flat group 
schemes over $\Spec R$ is a short exact sequence if its restriction to $\Spec
R\setminus\{\ideal{m}\}$ is a
short exact sequence. This is a local statement in the faithfully flat topology of
$\Spec R$ and thus to check it we can assume that $R=\hat{R}= W(k)[[T_1,T_2]/(p -
h)$ with $h\in (T_1,T_2)$ but $h \notin (p,T_1^p,T_2^p,T_1^{p-1}T_2^{p-1})$.
By Lemma  \ref{Le3} we can assume that $H_2$ and  $H_3$ are connected. Thus
$0\rightarrow H_1\rightarrow H_2\rightarrow H_3\rightarrow 0$ is a short exact
sequence, cf. Corollary \ref{C2}. \endproof

\section{Extending abelian schemes}

\begin{prop}\label{P7}
Let $R$ be a regular local ring of mixed characteristic $(0,p)$. 

\medskip
{\bf (a)} Let $U=\Spec R\setminus\{\ideal{m}\}$, where $\ideal{m}$ is the maximal
ideal of $R$. Let $\breve{A}$ be an abelian scheme over $U$. If the $p$-divisible
group of $\breve{A}$ extends to a $p$-divisible group over $\Spec R$, then
$\breve{A}$ extends uniquely to an abelian scheme $A$ over $\Spec R$.

\smallskip
{\bf (b)} If $R$ is
$p$-quasi-healthy, then $R$ is quasi-healthy.
\end{prop}
{\bf Proof.} The uniqueness part of (a) is well known (cf. \cite{R1}, Chapter IX,
Corollary 1.4). Part (a) is a particular case of either proof of \cite{V2},
Proposition 4.1 (see remark that starts the proof) or  \cite{V2}, Remark 4.2. Part
(b) follows from (a).\endproof

\begin{lemma}\label{Le6}
Let $S \twoheadrightarrow R$ be a ring epimorphism between local noetherian rings
whose kernel is an ideal
$\ideal{a}$ with
$\ideal{a}^2 = 0$ and $\depth_R \ideal{a} \geq 2$. We assume that $\depth R\geq 2$
and that $R$ is
quasi-healthy. Then $S$ is quasi-healthy as well. 
\end{lemma}
{\bf Proof.} Let $\ideal{m}$ and $\ideal{n}$ be the maximal ideals of $R$ and $S$
(respectively).
We set $U= \Spec R\setminus \{\ideal{m}\}$ and $V = \Spec S
\setminus \{\ideal{n}\}$.

Let
$\breve{B}$ be an abelian scheme over $V$ and let $\breve{A}$ be its
reduction over $U$. Then $\breve{A}$ extends uniquely to an abelian scheme $A$
over $\Spec R$. Let $\pi: A \twoheadrightarrow \Spec R$ be the projection. 
It is well-known that the set of liftings of $A$ with respect to $\Spec R
\hookrightarrow \Spec S$ is a trivial torsor under the group $H^0(\Spec R,
R^1\pi_* \underline{\Hom}(\Omega_{A/R}, \ideal{a}))$ and that the set
of liftings of $\breve{A}$ with respect to $V\twoheadrightarrow U$ is a trivial
torsor under the group 
$H^0(U, R^1\pi_* \underline{\Hom}(\Omega_{A/R}, \ideal{a}))$. As $\depth_R \ideal{a}
\geq 2$, the last two
groups are equal. Thus there exists a unique abelian scheme $B$ over $\Spec S$ which
lifts $A$ and
whose restriction to $V$ is $\breve{B}$. \endproof

\begin{prop}\label{P9}
Let $S$ be a complete noetherian local ring of mixed characteristic $(0,p)$. Let $S
\twoheadrightarrow R$ be a ring epimorphism with kernel
$\ideal{a}$ (thus $R$ is a complete noetherian local ring). We assume that there
exists a
sequence
of ideals 
\begin{displaymath}
   \ideal{a} = \ideal{a}_0 \supset \ideal{a}_1 \supset \ldots 
\end{displaymath}
such that the intersection of these ideals is $0$ and for all $i \geq 0$ we have
$\ideal{a}_i^2 \subset \ideal{a}_{i+1}$ and $ \depth_{S} \ideal{a}_i/
\ideal{a}_{i+1} \geq 2$. We also assume that $\depth R\geq 2$ and that $R$ is
quasi-healthy. Let $\ideal{n}$ be the maximal ideal of $S$.  
 Let $V :=\Spec S \setminus \{\ideal{n}\}$. 

\medskip
{\bf (a)}  Then each abelian scheme over $V$ that has a polarization extends to an abelian scheme over $\Spec S$.

\smallskip
{\bf (b)} We assume that $S$ is integral and geometrically unibranch (like $S$ is
normal). Then $S$ is quasi-healthy. 
\end{prop}
{\bf Proof.} By a well-known theorem we
have $S =
\lim\limits_{\leftarrow} S/\ideal{a}_i$. 

For (b) (resp. (a)) we have to show that each abelian scheme $\breve{A}$ (resp. each
abelian scheme $\breve{A}$ that has a polarization $\lambda_{\breve{A}}$) over $V$
extends to an abelian scheme over $\Spec S$. We write $U_{i} =
(\Spec S/\ideal{a}_i) \setminus \{\ideal{n}\}$. The topological space
underlying $U_i$ is independent of $i$ and it will be denoted by $U$.  
It is easy to see that $\depth S/\ideal{a}_i\geq 2$. By Lemma \ref{Le6} the ring
$S/\ideal{a}_i$ is quasi-healthy.
We denote by $\breve{A}_i$ the base change of $\breve{A}$ to $U_i$. Then
$\breve{A}_i$ extends uniquely to an abelian scheme $A_i$ over $\Spec
S/\ideal{a}_i$. If $S$ is integral and geometrically unibranch, then from \cite{R1},
Chapter XI, Theorem 1.4 we
get that $\breve{A}$ is projective over $V$. Thus from now on we can assume that
there exists a polarization
$\lambda_{\breve{A}}$ of $\breve{A}$ and we will not anymore differentiate between
parts (a) and (b). 

From the uniqueness of $A_i$ we easily get that the reduction
of $\lambda_{\breve{A}}$ modulo $\ideal{a}_i$ extends uniquely to a polarization
$\lambda_{A_i}$ of $A_i$ (this also follows from \cite{R1}, Chapter IX, Corollary
1.4). We get that the $A_i$'s inherit a compatible system of
polarizations. From this and the algebraization theorem of Grothendieck, we get that
there exists an abelian scheme $A$ over $\Spec S$ which lifts the $A_i$'s. 

Next we will prove that the $p$-divisible group $G$ of $A$ restricts
over $V$ to the $p$-divisible group $\breve{G}$ of $\breve{A}$. This implies that
$\breve{A}$ extends to an abelian scheme
over $\Spec S$ (cf. Proposition \ref{P7} (a)) which is then necessarily isomorphic
to $A$.

Let $G_i[m]$ be the kernel of the multiplication by $p^m : A_i
\rightarrow A_i$, and let $G[m]$ be the kernel of the multiplication 
by $p^m : A \rightarrow A$.

Let
$C_i[m]$ be the $S/\ideal{a}_i$-algebra of global functions on $G_i[m]$. Then $B[m]  = 
\lim\limits_{\leftarrow} C_i[m]$ is the $S$-algebra of global functions on $G[m]$.  

We write $\underline{\Spec}\;
 \breve{\mathcal{C}}[m] = \breve{G}[m]$, where $\breve{\mathcal{C}}[m]$
 is a finite $\mathcal{O}_{\breve{G}[m]}$-algebra. We have a natural homomorphism:
\begin{displaymath}
   H^0(V, \breve{\mathcal{C}}[m]) \rightarrow 
   H^0(U, (C_i[m])^{\sim}) = C_i[m].
\end{displaymath}
Here $(C_i[m])^{\sim}$ is the restriction to $U$ of the
$\mathcal{O}_{\Spec S/\ideal{a}_i}$-algebra associated to $C_i[m]$. The last
equality follows from the fact that $\depth S/\ideal{a}_i \geq 2$. This gives birth
to an
$S$-algebra homomorphism 
\begin{displaymath}
    H^0(V, \breve{\mathcal{C}}[m]) \rightarrow B[m].
\end{displaymath}
If we restrict it to a homomorphism between $\mathcal{O}_V$-algebras we obtain a
homomorphism of finite flat group schemes over $V$
\begin{displaymath}
   G[m]_{V} \rightarrow \breve{G}[m]
\end{displaymath}
and thus a homomorphism of corresponding $p$-divisible groups
$G_{V} \rightarrow \breve{G}$ over $V$. By construction this is an isomorphism if we
restrict it to $U$. As $V$ is connected, from the following proposition we conclude
that $G_{V} \rightarrow \breve{G}$ is an isomorphism. \endproof

\begin{prop}\label{P10}
Let $\beta: G \rightarrow G'$ be a homomorphism of $p$-divisible groups over a
noetherian scheme $X$. Then the set $Y$ of points $x \in X$ such that
$\beta_x$ is an isomorphism is open and closed in $X$. Moreover $\beta_Y$ is an
isomorphism.
\end{prop}
{\bf Proof.} It is clear that $\beta$ is an isomorphism if and only if
$\beta[1] : G[1] \rightarrow G'[1]$ is an isomorphism. Therefore 
the subfunctor of $X$ defined by the condition that $\beta_Y$ 
is an isomorphism is representable by an open subscheme $Z \subset X$.

By the theorems of Tate and de Jong on extensions of
homomorphisms between $p$-divisible groups, we get that the valuative criterion of
properness holds for $Z\rightarrow X$. Thus $Z$ is as well closed in $X$ and
therefore we can take $Y=Z$.
\endproof
\subsection{Counterexample for  the $p$-quasi-healthy context}\label{S4} 
Lemma \ref{Le6} does not hold for the $p$-quasi-healthy context
even in the simplest cases. Here is an elementary counterexample. We
take $R=W(k)[[T_1]]$. From Subsection \ref{S3} we get that $R$ is
$p$-quasi-healthy. We take $S=R[T_2]/(T_2^2)$. Let $V:=\Spec
S\setminus \{\ideal{n}\}$, where $\ideal{n}$ is the maximal ideal of
$S$. Let $\mathcal{O}:=R_{(T_1)}$; we have a
natural identification $V=\Spec S[\frac{1}{T_1}] \cup \Spec
\mathcal{O}[T_2]/(T_2^2)$. 

Let $\mathcal{E}_k$ be an elliptic curve over $\Spec k$. We can identify the
formal deformation space of $\mathcal{E}_k$ with $\Spf R$. Let $\mathcal{E}$ be the
elliptic curve over $\Spec R$ which is the  algebraization of the uni-versal
elliptic curve over $\Spf R$. Let
$\iota_1:R\hookrightarrow S[\frac{1}{T_1}]$ be the  
$W(k)$-monomorphism that maps $T_1$ to $T_1+T_2T_1^{-1}$; it lifts the natural
$W(k)$-monomorphism $R\hookrightarrow R[\frac{1}{T_1}]$.
The $W(k)$-monomorphism $\iota_1$ gives birth to an elliptic curve
$\mathcal{E}_1$ over $\Spec S[\frac{1}{T_1}]$ whose restriction
to $\Spec S[\frac{1}{T_1}]/(T_2)=\Spec R[\frac{1}{T_1}]$ extends to the
elliptic curve $\mathcal{E}$ over $\Spec R$.

We check that the assumption that $\mathcal{E}_1$ extends to an elliptic curve
$\mathcal{E}_V$ over $V$ leads to a contradiction. To $\mathcal{E}_V$ and
$\mathcal{E}_U$ correspond morphisms $U\rightarrow V\rightarrow \mathbb{A}^1$, where
$\mathbb{A}^1$
is the $j$-line over $\Spec W(k)$. As the topological spaces of $U$ and $V$ are
equal, we get that we have a natural factorization $U\rightarrow V\rightarrow \Spec
R\rightarrow
\mathbb{A}^1$, where we identify $\Spec R$ with the completion of $\mathbb{A}^1$ at
its $k$-valued point defined by $\mathcal{E}_k$. As $S=H^0(V,\mathcal{O}_V)$, the
image of $\iota_1$ is contained in $S$. Contradiction.

But the $p$-divisible group $\mathcal{G}_1$ of $\mathcal{E}_1$ extends to
a $p$-divisible group $\mathcal{G}_V$ over $V$. This is so as each
$p$-divisible group  
over $\Spec \mathcal{O}$ extends uniquely to an \'etale $p$-divisible group over 
$\Spec \mathcal{O}[T_2]/(T_2^2)$. 

Finally we check that the assumption that $\mathcal{G}_V$ extends to a $p$-divisible
group $\mathcal{G}_2$ over $\Spec S$ leads to a contradiction. Let $\mathcal{E}_2$
be the elliptic curve over $\Spec S$ which lifts $\mathcal{E}$ and whose
$p$-divisible group is $\mathcal{G}_2$. Let $\iota_2:R\rightarrow S$ be the
$W(k)$-homomorphism that defines $\mathcal{E}_2$. We check that the resulting two
$W(k)$-homomorphisms $\iota_1,\iota_2:R\rightarrow S[\frac{1}{T_1}]$ are equal. It
suffices
to show that their composites $\iota_3,\iota_4:R\rightarrow \widehat{S_{(p)}}$ with the
natural $W(k)$-monomorphism $S[\frac{1}{T_1}]\hookrightarrow
\widehat{S_{(p)}}=\widehat{R_{(p)}}[T_2]/(T_2^2)$ are equal (here $\widehat{\bigtriangleup}$ denotes the completion of the local ring $\bigtriangleup$). But this follows from
Serre--Tate deformation theory and the fact that the composites
$\iota_5,\iota_6:R\rightarrow \widehat{R_{(p)}}$ of $\iota_3,\iota_4$ with the
$W(k)$-epimorphism $\widehat{S_{(p)}}\twoheadrightarrow\widehat{R_{(p)}}$, are
equal. As the two $W(k)$-homomorphisms $\iota_1,\iota_2:R\rightarrow
S[\frac{1}{T_1}]$ are
equal, the image of $\iota_1$ is contained in $S$. Contradiction.

We conclude that $S$ is not $p$-quasi-healthy. A similar argument shows that for all
$n\ge 2$, the ring $R[T_2]/(T_2^n)$ is not $p$-quasi-healthy.

This counterexample disproves the claims on \cite{FC}, top of p. 184 on torsors of
liftings of $p$-divisible groups.

\subsection{Proofs of Theorem \ref{T1} and Corollary \ref{T2}}\label{S5}
Let $\ideal{m}$ be the maximal ideal of $R$. Let $U=\Spec R\setminus\{\ideal{m}\}$.
Let $d=\dim R$.

We prove Theorem \ref{T1}. We first assume that $d=2$. It is enough to show that $R$ is
$p$-quasi-healthy, cf. Proposition \ref{P7} (b). But this follows from
Subsection \ref{S3}. We next assume that $d\geq 3$.  We have to show that each
abelian scheme over $U$ extends uniquely to an abelian scheme over $\Spec R$. This
is a  local
statement in the faithfully flat topology of $\Spec R$ and thus to
check it we can assume that $R=\hat{R}$ is complete with 
algebraically closed residue class field $k$. We have an epimorphism
$R\twoheadrightarrow W(k)[[T_1,T_2]]/(p - h)$ where $h$ is a power series in the 
maximal ideal of $W(k)[[T_1,T_2]]$ whose reduction modulo
$(p,T_1^p,T_2^p,T_1^{p-1}T_2^{p-1})$ is non-zero. As $W(k)[[T_1,T_2]]/(p -
h)$ is $p$-quasi-healthy (cf. the case $d=2$), from Proposition \ref{P9} we get that
$R$ is quasi-healthy. This proves Theorem \ref{T1}. 

We prove Corollary \ref{T2}. Thus $\dim R = d \geq 2$ and $p
\notin \ideal{m}^p$. As in the previous paragraph we argue that we can assume that
$R=\hat{R}$ is complete with 
algebraically closed residue class field $k$. We write $R =
\mathfrak{S}/(p - h)$ where $h\in (T_1,\ldots,T_d)$ is such that its
reduction $\bar{h} \in \bar{\mathfrak{S}}=\mathfrak{S}/p\mathfrak{S}$
modulo $p$ is a power series of order $e = \ord(\bar{h})\leq p-1$. 
Due to Noether normalization theorem we can assume that $\bar{h}$ contains
the monom $T_1^e$. We set $R' := \mathfrak{S}/(p - h, T_3, \ldots,
T_d)$; if $d=2$, then $R'=R$.  From Proposition \ref{P9} applied to the epimorphism
$R \twoheadrightarrow R'$,
we get that it suffices to show that $R^\prime$ is quasi-healthy and
$p$-quasi-healthy. Thus in order not
to complicate the notations, we can assume that $d=2$ (i.e., $R=R^\prime$). As $e\le
p-1$, the reduction of $\bar h$ modulo the ideal $(T_1^p,T_2^p,T_1^{p-1}T_2^{p-1})$
is non-zero. Thus Corollary \ref{T2} follows from Theorem \ref{T1}.\endproof  

\subsection{Example}\label{S6} Let $R$ be a regular local of mixed characteristic
$(0,p)$ such that the strict completion of $R$ is isomorphic to
\begin{displaymath}
   C_k[[T_1,\ldots, T_d]]/ (p - T_1^{e_1}\cdot \ldots \cdot T_m^{e_m})
\end{displaymath}
where $C_k$ is a Cohen ring of the field $k$, where $1\le m\le d$, and where the
$m$-tuple $(e_1,\ldots,e_m)\in
\mathbb{N}^m$ has the
property that there exists a disjoint union
$\{1,\ldots,m\}=I_1\bigsqcup I_2$ for which we have $m_1:=\sum_{i\in
  I_1} e_i\in\{1,\ldots,p-1\}$ and $m_2:=\sum_{i\in I_2}
e_i\in\{0,\ldots,p-2\}$. 

To check that $R$ is 
quasi-healthy we can assume that the field $k$ is algebraically closed and that
$R=W(k)[[T_1,\ldots, T_d]]/ (p - T_1^{e_1}\cdot \ldots \cdot
T_m^{e_m})$ (thus $C_k=W(k)$). We consider the ring epimorphism $R\twoheadrightarrow
R^\prime:=W(k)[[T_1,T_2]]/(p-T_1^{m_1}T_2^{m_2})$ that maps  $T_i$ with $i\in I_1$
to $T_1$, that maps  $T_i$ with $i\in I_2$ to $T_2$, and that maps $T_i$ with $i>m$
to $0$. From Theorem \ref{T1} we get that $R$ is quasi-healthy. 

Concrete example: if $1\le m\le\min\{d,2p-3\}$ and if the strict completion of $R$ is
$C_k[[T_1,\ldots, T_d]]/ (p - T_1\cdot \ldots \cdot T_m)$, then $R$ is
quasi-healthy. From this and Corollary \ref{C1} we get that each \'etale scheme over
$\Spec O[T_1,\ldots, T_d]/ (p - T_1\cdot \ldots \cdot T_m)$ is healthy regular,
provided $O$ is a discrete valuation ring of mixed characteristic $(0,p)$ and index
of ramification $1$.

\subsection{Regular schemes which are not  ($p$-) healthy}\label{S7}
Let $R$ be a local regular ring of dimension $2$ and mixed characteristic $(0,p)$.
Let $\ideal{m}$ be the maximal ideal of $R$. Let $U=\Spec R\setminus\{\ideal{m}\}$.

The ring $R$ is ($p$-) quasi-healthy if and only if $\Spec R$ is ($p$-) healthy
regular. The
next lemma provides an easy criterion for when $R$ is not ($p$-) quasi-healthy.

\begin{lemma}\label{Le5}
We assume that there exists a homomorphism $H\rightarrow D$ of finite flat group
schemes
over $\Spec R$ which is not an epimorphism and whose restriction to $U$ is an
epimorphism.
Then $R$ is neither quasi-healthy nor $p$-quasi-healthy. 
\end{lemma} 
{\bf Proof.} We embed the Cartier dual $D^{\text{t}}$ of $D$ into an abelian scheme
$A$ over $\Spec R$, cf. \cite{BBM}, Theorem 3.1.1. The Cartier dual homomorphism
$H^{\text{t}}_{U} \rightarrow D^{\text{t}}_{U}$ is a
closed immersion. The abelian scheme $A_{U}/H^{\text{t}}_{U}$
over $U$ does not extend to an abelian scheme over $\Spec R$. Based on \cite{R1},
Chapter IX, Corollary 1.4, the argument for this is similar to the one used to prove
that (d) implies (c) in Lemma \ref{Le3}. Thus $R$ is not quasi-healthy. From
Proposition \ref{P7} (a) we get that the $p$-divisible group of
$A_{U}/H^{\text{t}}_{U}$ does not extend to $\Spec R$. Thus the fact that $R$ is not
$p$-quasi-healthy follows from either Lemma \ref{Le3} or Proposition \ref{P7} (a).
\endproof

\medskip
Based on Lemma \ref{Le5}, the following theorem adds many examples to the classical
example of Raynaud--Gabber. 

\begin{thm}\label{T4}
We consider the ring $\mathfrak{S} = W(k)[[T_1,T_2]]$. Let $h \in
(T_1,T_2)\setminus p\mathfrak{S}$. Let $R:=\mathfrak{S}/(p-h)$. Let
$\bar{\mathfrak{S}}:=\mathfrak{S}/p\mathfrak{S}$ and let $\bar{h}$ be the reduction
of $h$ modulo $p$.
We assume that one of the following three properties hold:

\medskip
{\bf (i)} The element $\bar h$ is divisible by $u^p$, where $u$ is a power series in
the maximal ideal
of $\bar{\mathfrak{S}}$ (the class of rings $R$ for which (i) holds includes
$O[[T]]$, where $O$ is a totally ramified
discrete valuation ring extension of $W(k)$ of index of ramification at least equal
to $p$).

\smallskip
{\bf (ii)} There exists a regular sequence $u,v$ in
$\bar{\mathfrak{S}}$ such that $u^{p-1}v^{p-1}$ divides $\bar{h}$.

\smallskip
{\bf (iii)} We can write $\bar h=(aT_1^p+bT_2^p+cT_1^{p-1}T_2^{p-1})c$, where
$a,b,c\in\bar{\mathfrak{S}}$.

\medskip
Then $R$ is neither $p$-quasi-healthy nor quasi-healthy.
\end{thm}
{\bf Proof.} It is enough to construct a homomorphism $\beta:H \rightarrow D$ of
connected finite flat group
schemes over $\Spec R$ which is not an epimorphism but whose
restriction to $U$ is an
epimorphism, cf. Lemma \ref{Le5}. 

We will construct $H$ and $D$ by specifying their Breuil modules $(M,\varphi)$ and
$(N, \tau)$ (respectively) associated to a standard frame for $R$. We first present
with full details the case (i) and then we will only mention what are the changes
required to be made for the other two cases (ii) and (iii). 

We assume that the condition (i) holds. We set $M = \bar{\mathfrak{S}}^3$ and we
identify $M^{(\sigma)}$ with
$\bar{\mathfrak{S}}^3$. We choose an element $t \in
\bar{\mathfrak{S}}$ such that $t$ and $u$ are a regular sequence in $\bar{\mathfrak{S}}$. 
We define the homomorphism $\varphi$ by the
following matrix:
 
\begin{displaymath}
\Gamma =  \left(
\begin{array}{rrr}
0 & 0 & u^p\\
t - t^pu^{p-1} & u & (u - t^{p-1})(t - t^pu^{p-1})\\
u^{p-1} & 0 & u^{p-1}(u - t^{p-1})
\end{array}
\right).
\end{displaymath}
It is easy to see that there exists a matrix $\Delta\in M_{3 \times
3}(\bar{\mathfrak{S}})$ such that 
\begin{displaymath}
   \Delta\Gamma = \Gamma\Delta = u^pI_3,
\end{displaymath}
where $I_3$ is the unit matrix. It follows that the cokernel of $\varphi$
is annihilated by $u^p$ and thus also by $\bar{h}$. Moreover the image of $\varphi$ is
contained in $(t,u)M$. Therefore
$(M,\varphi)$ is the Breuil module $H$ of a connected finite flat group
scheme over
$\Spec R$ annihilated by $p$. 

We set $N = \bar{\mathfrak{S}}$ and $N^{(\sigma)} = \bar{\mathfrak{S}}$
and we define $\tau$ as the multiplication by $u^p$. This defines
another connected finite flat group scheme $D$ over $\Spec R$ annihilated by
$p$.

One easily checks the following equation of matrices:
\begin{displaymath}
   (t^p,u^p,(tu)^p) \Gamma = u^p (t,u,tu).
\end{displaymath}
This equation shows that the $\mathfrak{S}$-linear map $M \rightarrow N$ defined
by the matrix $(t,u,tu)$ is a morphism of Breuil modules

\begin{displaymath}
   \alpha: (M,\varphi) \rightarrow (N,\tau).
\end{displaymath}
As $\alpha$ is not surjective, the homomorphism $\beta:H \rightarrow D$ associated to
$\alpha$ is not
an epimorphism. Let $\ideal{p} \ne \ideal{m}$ be a
prime ideal of $R$ which contains $p$. The base change
of $\alpha$ by  $\kappa_{\ideal{p}}:\mathfrak{S}
\rightarrow W(\kappa(\ideal{p})^{\perf})$ (of Section \ref{S0}) is an epimorphism as
the cokernel of
$\alpha $ is $k$. This implies that $\beta_U:H_U\rightarrow D_U$ is an
epimorphism.

We assume that the condition (ii) holds. The proof in this case is similar to the
case (i) but
with the definitions of $M$, $\Gamma$, $\tau$, and $M\rightarrow N$ modified as
follows. Let
$M=\bar{\mathfrak{S}}^2$. Let
\begin{displaymath}
\Gamma =  \left(
\begin{array}{rr}
u^{p-1} & 0\\
0 & v^{p-1}
\end{array}
\right).
\end{displaymath} 
Let $\tau$ be defined by $(uv)^{p-1}$. Let $M\rightarrow N$ be defined by the matrix
$(v\;
u)$.

We assume that the condition (iii) holds. Let
$M=\bar{\mathfrak{S}}^2$. Let
\begin{displaymath}
\Gamma =  \left(
\begin{array}{rr}
aT_1+cT_2^{p-1} & a T_2 \\
bT_1 & bT_2+cT_1^{p-1} 
\end{array}
\right).
\end{displaymath} 
Let $\tau$ be defined by $aT_1^p+bT_2^p+cT_1^{p-1}T_2^{p-1}$. Let $M\rightarrow N$
be defined by the matrix $(T_1\;T_2)$. The determinant of $\Gamma$ is $\bar h$. 
\endproof

\section{Integral models and N\'eron models} 

Let $O$ be a discrete valuation ring of mixed characteristic $(0,p)$ and of index of
ramification at most $p-1$. Let $K$ be the field of fractions of $O$. A flat
$O$-scheme $\bigstar$ is said to have the {\it extension property}, if for each 
$O$-scheme $X$ which is a  healthy regular scheme, every morphism $X_K\rightarrow
\bigstar_K$ of $K$-schemes extends uniquely to a morphism $X\rightarrow\bigstar$ of
$O$-schemes (cf. \cite{V1}, Definition 3.2.3 3)).

\begin{lemma}\label{Le4}
Let $Z_K$ be a regular scheme which is formally smooth over $\Spec K$. Then
there exists at most one regular scheme which is a formally smooth
$O$-scheme, which has the extension property, and whose fibre over $\Spec K$
is $Z_K$. 
\end{lemma}
{\bf Proof.} Let $Z_1$ and $Z_2$ be two regular schemes which are
formally smooth over $\Spec O$, which satisfy the identity
$Z_{1,K}=Z_{2,K}=Z_K$, and which have the extension property. Both
$Z_1$ and $Z_2$ are healthy regular schemes, cf. Corollary
\ref{C1}. Thus the identity $Z_{1,K}=Z_{2,K}$ extends naturally to
morphisms $Z_1\rightarrow Z_2$ and $Z_2\rightarrow Z_1$, cf. the fact that both $Z_1$
and $Z_2$ have the extension property. Due to the uniqueness part of
the extension property, the composite morphisms $Z_1\rightarrow Z_2\rightarrow Z_1$
and $Z_2\rightarrow Z_1\rightarrow Z_2$ are identity automorphisms. Thus the identity
$Z_{1,K}=Z_{2,K}$ extends uniquely to an isomorphism $Z_1\rightarrow
Z_2$.\endproof 
 
\begin{cor}\label{C5}
The {\it integral canonical
  models} of Shimura varieties defined in \cite{V1}, Definition 3.2.3
6) are unique, provided they are over the spectrum of a discrete valuation ring $O$ as
above. 
\end{cor}

Let $d\ge 1$ and $n\geq 3$ be natural numbers. We assume that $n$ is prime to $p$.
Let $\mathcal A_{d,1,n}$ be the Mumford moduli scheme over
$\Spec \mathbb{Z}[\frac{1}{n}]$ that parameterizes principally polarized
abelian scheme over $\Spec \mathbb{Z}[\frac{1}{n}]$-schemes which are of relative
dimension $d$ and
which have level-$n$ symplectic similitude structures (cf. \cite{MFK},
Theorems 7.9 and 7.10). For N\'eron models over Dedekind domains we
refer to \cite{BLR}, Chapter I, Subsection 1.2, Definition 1. 

\begin{thm}\label{T5}
 Let $\mathcal{D}$ be a Dedekind domain which is a flat
 $\mathbb{Z}[\frac{1}{n}]$-algebra. Let $\mathcal{K}$ be the field of fractions
 of $\mathcal{D}$. We assume that the following two things hold: 

\medskip\noindent
{\bf (i)} the only local ring of $\mathcal{D}$ whose residue class field has
characteristic $0$, is
$\mathcal{K}$;

\smallskip\noindent
{\bf (ii)} if $v$ is a prime of $\mathcal{D}$ whose residue class field has a prime
characteristic $p_v\in \mathbb{N}^*$, then the index of ramification of
the local ring of $v$ is at most $p_v-1$. 

\medskip
 Let $\mathcal N$ be a finite $\mathcal A_{d,1,n,\mathcal{D}}$-scheme which is a
projective,
smooth $\mathcal{D}$-scheme. Then $\mathcal N$ is the N\'eron model over
$\mathcal{D}$ of its generic
fibre $\mathcal N_{\mathcal{K}}$.
\end{thm}

{\bf Proof.}
Let $Y$ be a smooth $\mathcal{D}$-scheme. Let
$\delta_{\mathcal{K}}:Y_{\mathcal{K}}\rightarrow \mathcal N_{\mathcal{K}}$ be a
morphism of
${\mathcal{K}}$-scheme. Let $V$ be an open subscheme of $Y$ which contains
$Y_{\mathcal{K}}$ and all generic
points of fibres of $Y$ in positive characteristic and for which
$\delta_{\mathcal{K}}$ extends
uniquely to a morphism $\delta_V:V\rightarrow\mathcal N$ (cf. the projectiveness of
$\mathcal
N$). Let $(\mathcal B_V,\lambda_V)$ be the pull back to $V$ of the universal
principally polarized abelian scheme over  $\mathcal A_{d,1,n,\mathcal{D}}$ via the
composite
morphism $\nu_V:V\rightarrow\mathcal N\rightarrow\mathcal A_{d,1,n,\mathcal{D}}$.
From Corollary 
\ref{C1}  we
get that $\mathcal B_V$ extends uniquely to an abelian scheme $\mathcal B$ over $Y$.
From \cite{R1}, Chapter IX, Corollary 1.4 we get that $\lambda_V$ extends (uniquely)
to a polarization
$\lambda$ of $\mathcal B$. The level-$n$ symplectic similitude structure of 
$(\mathcal B_V,\lambda_V)$ defined naturally by $\nu_V$ extends uniquely to a
level-$n$ symplectic similitude structure of  $(\mathcal B,\lambda)$, cf. the
classical Nagata--Zariski purity theorem. Thus $\nu_V$
extends uniquely to a morphism
$\nu:Y\rightarrow\mathcal A_{d,1,n,\mathcal{D}}$. As $Y$ is a normal scheme, as
$\mathcal N$
is finite
over $\mathcal A_{d,1,n,\mathcal{D}}$,  and as $\nu$ restricted to $V$ factors
through
$\mathcal N$, the morphism $\nu$ factors uniquely through a morphism
$\delta:Y\rightarrow\mathcal
N$ which extends $\delta_V$ and thus also $\delta_{\mathcal{K}}$. Hence the Theorem
follows from the very definition of N\'eron models.\endproof

\medskip\medskip
\noindent
{\bf Remarks.} (a) From Theorem  \ref{T5} and \cite{V3}, Remark 4.4.2 and Example 4.5
we get that there exist plenty of N\'eron models over $O$ whose generic fibres are not
finite schemes over torsors of smooth
schemes over $\Spec K$.

We can take $\mathcal N$ to be the pull back to $O$ of those N\'eron models
of Theorem \ref{T5}
 whose generic fibres have the above property (cf. \cite{V3}, Remark 4.5). If $p>2$
and $e=p-1$, then these N\'eron models
$\mathcal N$ are new (i.e., their existence does not follow
from \cite{N}, \cite{BLR}, \cite{V1}, \cite{V2}, or \cite{V3}). 

(b) One can use Theorem \ref{T4} (i) and Artin's approximation theorem to show
that Theorem \ref{T5} does not hold in general if there exists a prime $v$ of
$\mathcal{D}$
whose residue class field has a prime characteristic $p_v\in \mathbb{N}^*$ and
whose index of ramification is at least $p_v$. Counterexamples
can be obtained using integral models of projective Shimura
varieties of PEL type, cf. \cite{V3}, Corollary 4.3.

\bigskip

\hbox{Adrian Vasiu,\;\;\;E-mail: adrian@math.binghamton.edu}
\hbox{Address: P.O. Box 6000, Department of Mathematical Sciences, Binghamton University,}
\hbox{Binghamton, New York 13902-6000, U.S.A.}

\hbox{}
\hbox{Thomas Zink,\;\;\;E-mail: zink@math.uni-bielefeld.de}
\hbox{Address: Fakult\"at f\"ur Mathematik, Universit\"at Bielefeld,} \hbox{P.O. Box
100 131, D-33 501 Bielefeld, Germany.}

\end{document}